\documentclass[12pt]{article}
\usepackage{pifont}
\usepackage{psfrag,float,subfigure}
\usepackage{float,amsmath,amssymb,amsthm,mathrsfs,graphicx,subfigure,mathrsfs}
\usepackage{caption,enumerate}
\usepackage[numbers,sort&compress]{natbib}
\usepackage[colorlinks]{hyperref}
\usepackage[scale=0.8,a4paper]{geometry}

\newtheorem{theorem}{Theorem}[section]
\newtheorem{lemma}{Lemma}[section]

\newtheorem{corollary}[lemma]{Corollary}

\newtheorem{prop}[lemma]{Proposition}

\newtheorem{observation}[lemma]{Observation}

\title{Removable edges in near-bipartite bricks}
\author{{Yipei Zhang$^1$, Fuliang Lu$^2$\thanks{Corresponding author. \newline\indent\indent\!\!\! Email addresses: zyipei@163.com, flianglu@163.com, wangxiumei@zzu.edu.cn, yuanjj@zzu.edu.cn.}, Xiumei Wang$^1$, Jinjiang Yuan$^1$}\\
{\small $^{1}$  School of Mathematics and Statistics, Zhengzhou University,}
{\small Zhengzhou 450001, China}\\
{\small $^{2}$ School of Mathematics and Statistics,
Minnan Normal University, Zhangzhou 363000, China}
}
\date{}\makeatother
\begin{document}
\maketitle

\begin{abstract}
An edge $e$ of a matching covered graph $G$ is \emph{removable} if $G-e$ is also matching covered. The notion of removable edge arises in connection with ear decompositions of
matching covered graphs introduced by Lov\'asz and Plummer. A nonbipartite matching covered graph $G$ is a \emph{brick} if it is free of nontrivial tight cuts. Carvalho, Lucchesi, and Murty proved that every brick other than $K_4$ and $\overline{C_6}$ has at least $\Delta-2$ removable edges.
A brick $G$ is \emph{near-bipartite} if it has a pair of edges $\{e_1,e_2\}$ such that $G-\{e_1,e_2\}$ is a bipartite matching covered graph.
In this paper, we show that in a near-bipartite brick $G$ with
at least six vertices, every vertex of $G$, except at most six vertices of degree three contained in two disjoint triangles, is incident with at most two nonremovable edges;
consequently, $G$ has at least $\frac{|V(G)|-6}{2}$ removable edges.
Moreover, all graphs attaining this lower bound are characterized.
\end{abstract}


\section{Introduction}

Graphs considered in this paper may have multiple edges, but no
loops. We follow \cite{BM08} for undefined notation and terminology.
A connected graph $G$ is \emph{$k$-extendable} if each set of $k$ independent edges extends to a perfect matching of $G$.
An edge $e$ of a graph $G$ is \emph{admissible} if $G$ has a perfect matching that contains $e$, and \emph{nonadmissible} otherwise.
A connected nontrivial graph  is \emph{matching covered} if each of its edges is admissible. Matching covered graphs are also called 1-extendable \cite{LP86}.
For $X\subseteq V(G)$, by $\partial(X)$ we mean the {\it edge cut} of $G$, which is the set of edges of $G$ with one end in $x$ and the other in  $\overline{X}$, where $\overline{X}=V(G)\backslash X$; by $G/X\rightarrow x$ or simply $G/X$ we mean the graph obtained by contracting $X$ to a single vertex $x$, the graph $G/\overline{X}\rightarrow \overline{x}$ or simply $G/\overline{X}$ is defined analogously.

Let $G$ be a matching covered graph.
An edge cut $C=\partial(X)$ of  $G$ is {\it tight} if $|M\cap C|=1$ for each perfect matching $M$ of $G$ and is {\it separating} if $G/X$ and $G/{\overline X}$ are matching covered.
A matching covered graph that is free of nontrivial tight cuts is a {\it brace} if
it is bipartite, and a \emph{brick} if it is nonbipartite. Edmonds et al. \cite{ELP82} showed that a graph $G$ is a \emph{brick} if and only if $G$ is 3-connected and for any two distinct vertices $x$ and $y$ of $G$, $G-x-y$ has a perfect matching. A brick is \emph{solid} if it is free of nontrivial separating cuts.  We denote the number of vertices of $G$ by $n$.
An edge $e$ of $G$ is \emph{removable} if $G-e$ is also matching covered.
A pair of edges $\{e_1,e_2\}$ is {\it a removable doubleton} of $G$ if neither $e_1$ nor $e_2$ is removable in $G$ but $G-\{e_1,e_2\}$ is matching covered. The notion of removable edge arises in connection with ear decompositions of matching covered graphs introduced by Lov\'asz and Plummer.
The existence of removable edges, especially of special types, plays an important role in the generation of some bricks, see \cite{NK19,NK20}.
Lov\'asz \cite{Lovasz87} first showed the existence of removable edges of bricks other than $K_4$ and $\overline{C_6}$.
Carvalho, Lucchesi and Murty \cite{CLM99, CLM12} proved that every brick other than $K_4$ and $\overline{C_6}$ has at least $\Delta-2$ removable edges; in every solid brick $G$ with six or more  vertices, each vertex is incident with  at most two nonremovable edges, consequently, $G$ has at least $\frac{n}{2}$ removable edges.
Zhai, Lucchesi and Guo \cite{zhai} showed that every matching covered graph $G$
has at least $m(G)$   removable classes (including removable edges and removable doubletons), where $m(G)$ denotes the minimum number of perfect matchings needed to cover all edges of $G$.
For bipartite matching covered graphs, He et al. \cite{HWYZ2019} gave a characterization of graphs each of whose edges is removable.

\begin{figure}[h]
 \centering
 \includegraphics[width=0.85\textwidth]{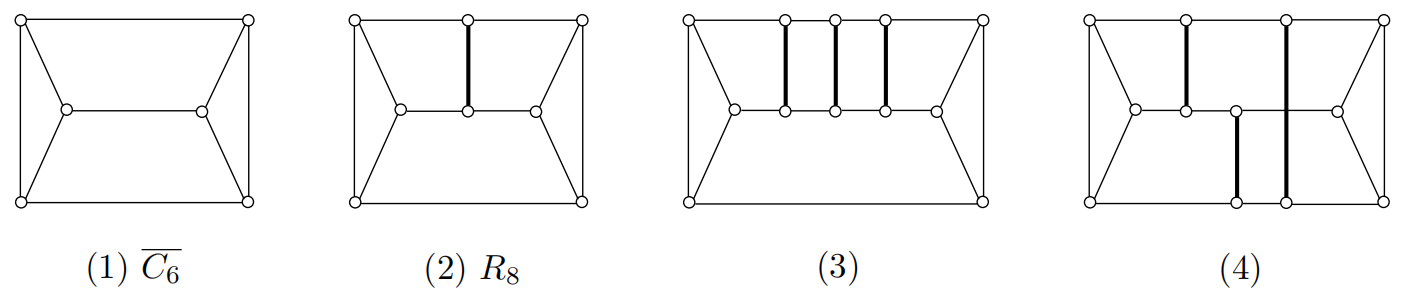}\\
 \caption{The four near-bipartite tri-ladders}\label{fig1}
\end{figure}

A nonbipartite matching covered graph $G$ is {\it near-bipartite} if it has a pair of edges $\{e_1,e_2\}$ such that $G-\{e_1,e_2\}$ is a bipartite matching covered graph; such a pair $\{e_1,e_2\}$ is referred to as a \emph{removable doubleton} of $G$.
The graphs $K_4$ and $\overline{C_6}$ are the only simple near-bipartite bricks on four and six vertices, respectively, each of which has three removable doubletons but no removable edges. The significance of near-bipartite graphs arises from the theory of ear decompositions. Fischer and Little  characterized Pfaffian near-bipartite graphs in \cite{FL}. Kothari \cite{NK19}, and Kothari and  Carvalho \cite{NK20} investigated generation procedures
for near-bipartite bricks  and  simple near-bipartite bricks, respectively.
Inspiring by the structure with respect to nonremovable edges in solid bricks,  we consider near-bipartite bricks in this paper. The main results are stated as follows.

\begin{theorem}\label{th:main1}
Let $G$ be a near-bipartite brick other than $K_4$. Then every vertex of $G$,
except at most six vertices of degree three contained in two disjoint triangles, is incident with at most two nonremovable edges.
\end{theorem}

\begin{theorem}\label{th:main2}
Every near-bipartite brick $G$ other than $K_4$ has at least $\frac{n-6}{2}$ removable edges, and near-bipartite tri-ladders (see Section 4) are the only graphs attaining the lower bound.
\end{theorem}


This paper is organized as follows. In Section 2, we present some basic results. In Section 3 and 4,  we give proofs of Theorems \ref{th:main1} and \ref{th:main2}, respectively.

\section{Preliminaries}
We begin with some notations. Let $G$ be a graph with the vertex set $V(G)$ and the edge set $E(G)$. By $V(e)$ we mean the set of the two ends of the edge $e$. For $X\subseteq V(G)$,  by $N_G(X)$, or simply $N(X)$, we mean the set of vertices that are not in $X$ but have neighbours in $X$; by $G[X]$ we mean the subgraph of $G$ induced by $X$.
For $X,Y\subseteq V(G)$, by $E_G[X,Y]$, or simply $E[X,Y]$, we mean  the set of edges of $G$ with one end in $X$ and the other end in $Y$.
If $Y=X$, we use $E(X)$ instead of $E[X,X]$.
Then the edge cut $\partial(X)=E[X,\overline{X}]$.
An edge cut $\partial(X)$ is a {\it $k$-cut} if $|\partial(X)|=k$; is {\it trivial} if either $|X|=1$ or $|\overline{X}|=1$, and {\it nontrivial} otherwise.

\subsection{Basic results}
Tutte \cite{Tutte47} proved that  a graph $G$ has a perfect matching  if and only if $o(G-S)\leq|S|$ for every $S\subseteq V(G)$, where $o(G-S)$ denotes  the number of odd components of $G-S$. A nonempty subset $S$ of $V(G)$ is a {\it barrier} of $G$ that has a perfect matching if $o(G-S)=|S|$. A graph $G$ is {\it factor-critical} if, for each vertex $v$ of $G$,  $G-v$ has a perfect matching.  Using Tutte's theorem, Lemma \ref{lem-barrier1} may be easily derived (see \cite{CLM12}).

\begin{lemma}\label{lem-barrier1}
Let $G$ be a graph with a perfect matching. Then
\vspace{-8pt}
\begin{enumerate}[(\romannumeral1)]
\setlength{\itemsep}{-1ex}
\item if $G$ is  a brick, then every barrier of $G$ is trivial,
\item an edge $e$ of $G$ is admissible if and only if $G$ has no barriers containing $V(e)$,  and
\item  for each  maximal barrier $B$ of $G$, all components of $G-B$ are factor-critical.
\end{enumerate}
\end{lemma}

\begin{lemma}[\cite{CLM12}]\label{lem-barrier-cap}
Let $G$ be a brick, and let $f_1$ and $f_2$ be two adjacent edges of $G$. If, for $i=1,2$, $S_i$ is a barrier of $G-f_i$, then $|S_1\cap S_2|\leq 1$.
\end{lemma}

\begin{lemma}\label{lem-barrier2}
Let $G$ be a  graph with a perfect matching, let $S$ and $S'$ be two subsets of $V(G)$ such that $N(S)\subseteq S'$, $S\cap S'=\emptyset$ and $|S'|\leq |S|+1$. If $S$ is an independent set of $G$, then $S'$ is a  barrier of $G$.
\end{lemma}

\begin{proof}
Since $N(S)\subseteq S'$, $S\cap S'=\emptyset$ and $S$ is an independent set, we have $o(G-S')\geq |S|$. Furthermore, since $G$ has a perfect matching, by Tutte's theorem and the assumption that $|S'|\leq |S|+1$, we have
$$|S'|-1\leq|S|\leq o(G-S')\leq |S'|,$$
Thus $o(G-S')=|S'|$ or $|S'|-1$.
Since $|V(G)|$ is even, 
$o(G-S')$ and  $|S'|$ have the same parity.  Thus $o(G-S')=|S'|$ and then $S'$ is a  barrier of $G$.
\end{proof}

Note that it is possible that $S'=N(S)$  in Lemma \ref{lem-barrier2}, i.e., if $S$ is an independent set of a graph $G$ with a perfect matching and $|N(S)|\leq |S|+1$,  then $N(S)$ is a  barrier of $G$.

\begin{lemma}[\cite{LP86}]\label{lem-MC-N(S)}
Let $H[U,W]$ be a bipartite graph with four or more vertices, where
$|U|=|W|$. Then $H$ is matching covered if and only if $|N(S)|\geq |S|+1$ for every nonempty proper subset $S$ of $U$.
\end{lemma}

\begin{lemma}[\cite{LovaszVem}]\label{MC-nonremovable}
Let $H[U,W]$ be a bipartite matching covered  graph, and $u$ a
vertex of $G$ of degree $d$, where $u\in U$  and $d\geq3$. If $uw_1,uw_2,\ldots,uw_r$, $0<r\leq d$, are nonremovable edges of $H$ incident with $u$, then there exist partitions $(U_0,U_1,\ldots,U_r)$ of $U$ and $(W_0,W_1,\ldots,W_r$) of $W$ such that $u\in U_0$ and, for $i\in \{1,2,\cdots,r\}$,
(a) $|U_i|=|W_i|$, (b) $w_i\in W_i$ and (c) $N(W_i)=U_i\cup \{u\}$; in particular, $uw_i$ is the only edge in $E[U_0,W_i]$.
\end{lemma}

From Lemma \ref{MC-nonremovable}, we can obtain the following lemma.

\begin{lemma}\label{lem-4-cycle}
Let $H$ be a bipartite matching covered  graph, and $u$ a vertex of $H$ with
degree three or more. If $f_1$ and $f_2$ are two edges incident with $u$ that lie in a 4-cycle,  then at least one of $f_1$ and $f_2$ is removable.
\end{lemma}

\subsection{Two  types of nonremovable edges}
For the rest of Section 2 and Section 3, we assume that
\vspace{-9pt}
\begin{enumerate}[(\romannumeral1)]
\setlength{\itemsep}{-1ex}
\item $G$ is a near-bipartite brick with a removable doubleton $\{e_1,e_2\}$, and
\item $(U,W)$ is the bipartition of $H$ such that $e_1$ connects two vertices of $U$ and $e_2$ connects two vertices of $W$, where $H=G-\{e_1,e_2\}$.
\end{enumerate}
\vspace{-6pt}
Note that both $U$ and $W$ are stable set of $H$, $e_1$ is the only edge in $E_G(U)$ and $e_2$ is the only edge in $E_G(W)$.
Unless otherwise specified, we use $N(X)$, $E(X)$, and $E[X,Y]$ for $N_G(X)$, $E_G(X)$, and $E_G[X,Y]$, respectively, where $X,Y\subseteq V(G)$.
Then $E(U)=\{e_1\}$ and $E(W)=\{e_2\}$.

\begin{lemma}\label{lem-near-N(S)}
Let $S$ be a subset  of $U$ (or $W$) that contains at most one end of $e_1$ (or $e_2$). If $N(S)\subseteq S'$ and $|S'|\geq 2$, then  $|S'|\geq|S|+2$.
\end{lemma}

\begin{proof}
If $S=\emptyset$, the assertion is trivial. Now suppose that $S\neq\emptyset$. Since $H$ is a bipartite matching covered  graph, Lemma \ref{lem-MC-N(S)} implies that $|N_H(S)|\geq |S|+1$. If $S$  contains exactly one end of $e_1$, then $|N(S)|=|N_H(S)|+1\geq |S|+2$. If $S$ contains no ends of $e_1$, then $|N(S)|=|N_H(S)|\geq |S|+1$.  When $|N(S)|=|S|+1$, Lemma \ref{lem-barrier2} implies that $N(S)$ is a nontrivial barrier of $G$. Since $G$ is a brick, Lemma \ref{lem-barrier1}(i) implies that $G$ has no nontrivial barriers, a contradiction. Thus $|N(S)|\geq |S|+2$. Since $N(S)\subseteq S'$, we have $|S'|\geq|N(S)|\geq |S|+2$.
\end{proof}

Since $\{e_1,e_2\}$ is a removable doubleton of $G$,  none of $e_1$ and $e_2$ is removable in $G$.
Let $e$ be a nonremovable edge of $G$ such that $e\notin \{e_1,e_2\}$.
We say that $e$  is of {\it type I} if $e$ is removable in $H$, and  of {\it type II} otherwise.


\begin{lemma}\label{leme1e2e}
If $e$ is  an edge  of type I,  then
\vspace{-8pt}
\begin{enumerate}[(\romannumeral1)]
\setlength{\itemsep}{-1ex}
\item (\cite{CLM99}) there exist partitions $(A_1,A_2)$ of $U$ and $(B_1,B_2)$ of $W$ such that $|B_1|=|A_1|+1$, $|A_2|=|B_2|+1$, $\{e_1\}=E(A_2)$, $\{e_2\}= E(B_1)$ and $\{e\}=E[A_1,B_2]$, and
\item $B_1$ and $A_2$ are barriers of $G-e$.
\end{enumerate}
\end{lemma}

\begin{proof} (ii) If $e$ is of type I, by (i),  $N_{G-e}(A_1)\subseteq B_1$. Since $|B_1|=|A_1|+1$, Lemma \ref{lem-barrier2} implies that $B_1$ is a barrier of $G-e$. Similarly, $A_2$ is also a barrier of $G-e$.
\end{proof}

Now assume that $e$ is of type II, by the definition, $e$ is nonremovable in both $G$ and $H$.
So both $G-e$ and $H-e$ have nonadmissible edges. Note that
if $h$ is a nonadmissible edge of $H-e$ and   is admissible in $G-e$, then the perfect matching of $G-e$ containing $h$ contains  $e_1$ and $e_2$.
Therefore, there exists an edge $e^*$ that is nonadmissible in both $H-e$  and $G-e$. Otherwise,  $G-e$ is matching covered, a contradiction.
By Lemma \ref{lem-barrier1}(ii), $G-e$ has a barrier containing $V(e^*)$.
Let $B$ be  a maximal such barrier. 
Lemma \ref{lem-barrier1}(iii) implies that each component of $G-e-B$ is factor-critical. In particular, each bipartite component is trivial, so $G-e-B$ has at most two nontrivial components, one contains $e_1$ and the other contains  $e_2$.

Let $U_2$  and  $W_1$ be the sets of the vertices in trivial components of $G-e-B$ that lie in $U$ and $W$, respectively.
Then $U_2$  and  $W_1$  are independent set of $G$.
Let $u$ and $w$ be the two ends of $e$ such that $u\in U$ and $w\in W$.
Since $e^*$ is admissible in $G$, $e$ connects two distinct components of  $G-e-B$, i.e.,  $u$ and $w$ lie in distinct components of $G-e-B$.
Let $U_1=B\cap U$ and $W_2=B\cap W$. Then $u\notin  U_1, w\notin W_2,$
\begin{equation}
e^*\in E[U_1,W_2], N_{H-e}(W_1)\subseteq U_1 \mbox{ and } N_{H-e}(U_2)\subseteq W_2.\label{e*}
\end{equation}

Let $\omega$ be the number of nontrivial components of $G-e-B$. Then $\omega\leq 2$ and
\begin{equation}
|B|=|U_1|+|W_2|=|W_1|+|U_2|+\omega.\label{Bw}
\end{equation}
Assume that $G_1$ and $G_2$ are the two nontrivial components of  $G-e-B$. When $\omega=0$, both $G_1$ and $G_2$ are null. When  $\omega=1$, for convenience, we assume that  $G_1$ is null if $|U_1|=|W_1|$ and $G_2$ is null otherwise.
Let $U_{i+2}=V(G_i)\cap U$ and $W_{i+2}=V(G_i)\cap W$, $i=1,2$.
Then $(U_1,U_2,U_3,U_4)$ is a partition of $U$ and $(W_1,W_2,W_3,W_4)$ is a partition of $W$.
Note that \begin{equation}N_{H-e}(U_i)\subseteq W_2\cup W_i \mbox{ and }N_{H-e}(W_i)\subseteq U_1\cup U_i \mbox{ for } i=3,4.\label{UW34}\end{equation}

\begin{observation}\label{ob}
When $\omega=0$, $U_3=W_3=U_4=W_4=\emptyset$;
when  $\omega=1$,  $U_3=W_3=\emptyset$ if $|U_1|=|W_1|$ and $U_4=W_4=\emptyset$ otherwise.
\end{observation}



\subsection{Properties}

\begin{prop}\label{prop} $U_1\neq\emptyset$, $ W_2\neq \emptyset$,  $|U_1|\geq |W_1|$ and $|W_2|\geq |U_2|$. \end{prop}

\begin{proof}
Since $ e^*\in E[U_1,W_2]$, we have $U_1\neq\emptyset$ and $ W_2\neq \emptyset$. To show $|U_1|\geq |W_1|$ and $|W_2|\geq |U_2|$, it suffices to consider the case $W_1\neq \emptyset$ and $U_2\neq \emptyset$. Note that $e$ is incident with at most one vertex of $W_1$ and $U_2$, respectively.
Since $H$ is a bipartite matching covered  graph, Lemma \ref{lem-MC-N(S)}  and (\ref{e*}) implies that $|U_1|+1\geq |N_H(W_1)|\geq |W_1|+1$ and $|W_2|+1\geq |N_H(U_2)|\geq |U_2|+1$, so $|U_1|\geq |W_1|$ and $|W_2|\geq |U_2|$.
\end{proof}

\begin{prop}\label{prop0} If $\omega=0$, then (i) $W_1\neq \emptyset$ and $U_2\neq \emptyset$; (ii) $|U_1|=|W_1|$ and $|U_2|=|W_2|$; (iii) $\{e\}= E[W_1,U_2]$,  $\{e_1\}=E[U_1,U_2]$ and $\{e_2\}=E[W_1,W_2]$.
\end{prop}

\begin{proof}
Since $\omega=0$, we have $u\in U_2$ and $w\in W_1$, so $U_2\neq \emptyset$ and $W_1\neq \emptyset$. Since $e$ is the only edge connecting two distinct components of $G-e-B$, $\{e\}=E[U_2, W_1]$ and $V(e_i)\cap B\neq \emptyset$ for $i=1,2$.  By (\ref{Bw}),  $|U_1|+|W_2|=|W_1|+|U_2|$.
By Proposition \ref{prop},  $|U_1|=|W_1|$ and $|U_2|=|W_2|$. Then $U_1\neq\emptyset$.
If $W_1$ contains no ends of $e_2$, then $N(W_1)\subseteq U_1\cup \{u\}$. By Lemma \ref{lem-near-N(S)}, $|U_1|\geq|W_1|+1$, a contradiction. Thus $\{e_2\}=E[W_1,W_2]$. Similarly, $\{e_1\}=E[U_1,U_2]$.
\end{proof}

\begin{prop}\label{prop1} Assume that $\omega=1$. Then
\vspace{-8pt}
\begin{enumerate}[(\romannumeral1)]
\setlength{\itemsep}{-1ex}
\item $|U_1|=|W_1|+1$ and $|U_2|=|W_2|$, or  $|U_1|=|W_1|$ and $|W_2|=|U_2|+1$;
\item if $|U_1|\neq|W_1|$, then $U_2\neq\emptyset$,  $\{e_1\}=E[U_1,U_2]$, $\{e_2\}=E(W_3)$, $|W_3|=|U_3|+1$, and $\{e\}= E[U_2, W_1]$ if $W_1\neq \emptyset$ and $\{e\}= E[U_2,W_3]$ otherwise;
\item if  $|U_1|=|W_1|$, then $W_1\neq\emptyset$,  $\{e_1\}=E(U_4)$, $\{e_2\}=E[W_1,W_2]$, $|U_4|=|W_4|+1$, and $\{e\}= E[W_1, U_2]$ if $U_2\neq\emptyset$ and  $\{e\}= E[W_1,U_4]$ otherwise.
\end{enumerate}
\end{prop}

\begin{proof}
 By (\ref{Bw}), $|U_1|+|W_2|=|W_1|+|U_2|+1$. (i) follows from Proposition \ref{prop}.

(ii) If $|U_1|\neq|W_1|$, by Observation \ref{ob} and (i), $U_4=W_4=\emptyset$ and $|U_2|=|W_2|$. By Proposition \ref{prop},  $W_2\neq\emptyset$, so $U_2\neq\emptyset$.
If at most one edge in $\{e,e_1\}$ has an end in $U_2$, say $e$, then $N(U_2)\subseteq W_2\cup \{w\}$. By Lemma \ref{lem-near-N(S)}, $|W_2|\geq|U_2|+1$, a contradiction.
Thus each edge of $\{e,e_1\}$ has exactly one end  in $U_2$. This implies that $\{e\}= E[U_2, W_1\cup W_3]$ and  $\{e_1\}=E[U_1,U_2]$, so  $\{e_2\}=E(W_3)$.
Recall that $G_1$ is factor-critical. Then $|W_3|=|U_3|+1$.
If $\{e\}=E[U_2,W_3]$, then $N(W_1)\subseteq U_1$. Since $|U_1|=|W_1|+1$, Lemma \ref{lem-barrier2} implies that $U_1$ is a barrier of $G$. By Lemma \ref{lem-barrier1}(i), $U_1$ is a singleton and $W_1=\emptyset$. (ii) holds.

If $|U_1|=|W_1|$, then $|U_2|\neq|W_2|$. Analogously, (iii) holds.
\end{proof}

Recall that when  $\omega=2$, $e_1$ and $e_2$ lie in two nontrivial components of $G-e-B$, respectively.  Assume, without loss of generality, that $e_1\in E(G_2)$ and $e_2\in E(G_1)$.
Recall that $\{e_1\}=E(U)$ and $\{e_2\}=E(W)$. Then
\begin{equation}\label{e12} \{e_1\}=E(U_4) \mbox{ and } \{e_2\}=E(W_3).
\end{equation}

\begin{prop}\label{prop2} Assume that $\omega=2$. Then
\vspace{-8pt}
\begin{enumerate}[(\romannumeral1)]
\setlength{\itemsep}{-1ex}
\item $|U_i|=|W_i|+1$ for $i=1, 4$ and  $|W_i|=|U_i|+1$ for $i=2,3$;
\item when $U_2\neq \emptyset$, $\{e\}=E[U_2,W_1]$ if $W_1\neq\emptyset$ and  $\{e\}= E[U_2,W_3]$ otherwise;
\item when $U_2= \emptyset$,  $\{e\}= E[U_4,W_1]$  if $W_1\neq\emptyset$ and $\{e\}= E[U_4,W_3]$ otherwise.
\end{enumerate}
\end{prop}

\begin{proof}
(i) By (\ref{e*}),  $N(W_1)\subseteq U_1\cup \{u\}$.  By Proposition \ref{prop}, $U_1\neq\emptyset$. By Lemma \ref{lem-near-N(S)}, $|U_1|\geq|W_1|+1$. Similarly, $|W_2|\geq|U_2|+1$. By (\ref{Bw}), $|U_1|+|W_2|=|W_1|+|U_2|+2$, so $|U_1|=|W_1|+1$ and $|W_2|=|U_2|+1$.
Recall that both $G_1$ and $G_2$ are factor-critical. By (\ref{e12}), we have $|U_4|=|W_4|+1$ and $|W_3|=|U_3|+1$.

To show (ii) and (iii), we first  claim that $V(e)\cap (U_3\cup W_4)=\emptyset$.
Let $B'= U_1\cup W_2\cup U_3\cup W_4$.
In fact, if $V(e)\cap (U_3\cup W_4)\neq\emptyset$,  then $\overline {B'}$ is an independent set of $H$.
By (i),  $|\overline{B'}|=|B'|$. Since $N_H(\overline {B'})\subseteq B'$,  by Lemma \ref{lem-barrier2}, $B'$ is a barrier of $H$.
By (\ref{e*}), $e^*\in E[U_1,W_2]$, so $e^*\in E(B')$.
Lemma \ref{lem-barrier1}(ii) implies that $e^*$ is nonadmissible in $H$, a contradiction. The claim holds.

(ii) When $U_2\neq \emptyset$, by (i), $|W_2|=|U_2|+1\geq2$.  If $u\notin U_2$, then $N(U_2)\subseteq W_2$. By Lemma \ref{lem-barrier2}, $W_2$ is a nontrivial barrier of $G$, a contradiction. Thus  $u\in U_2$. Likewise, if $W_1\neq\emptyset$, then $|U_1|\geq2$ and  $w\in W_1$. Thus $\{e\}=E[U_2,W_1]$.
If  $W_1=\emptyset$, then  $w\in W_3\cup W_4$. By the above claim,  $w\in W_3$, so $\{e\}= E[U_2,W_3]$.

(iii) When $U_2= \emptyset$, by the above claim again,  $u\in U_4$. If $W_1=\emptyset$,  then $\{e\}=E[U_4,W_3]$. If $W_1\neq \emptyset$,  by the same reason as the case $U_2\neq \emptyset$, then $w\in W_1$, so $\{e\}=E[U_4,W_1]$.   (iii) holds.
\end{proof}

\begin{corollary}\label{cor-P}
\begin{enumerate}[(\romannumeral1)]
\setlength{\itemsep}{-1ex}
\item If  $U_2\neq \emptyset$, then $u\in U_2$; otherwise,  $u\in U_4$ and $|W_2|=1$.
\item If $W_1\neq \emptyset$, then $w\in W_1$; otherwise, $w\in W_3$.
\item If $U_4\neq \emptyset$, then $\{e_1\}=E(U_4)$; otherwise, $\{e_1\}=E[U_1,U_2]$.
\item If $W_3\neq \emptyset$, then $|U_1|=|W_1|+1$ and $\{e_2\}=E(W_3)$; otherwise,  $|U_1|=|W_1|$, $\{e_2\}=E[W_1, W_2]$ and $U_3=\emptyset$.
\item $|U_2\cup U_4|=|W_2\cup W_4|$.
\end{enumerate}
\end{corollary}

\begin{proof}
Since $e=uv$ with $u\in U$ and $w\in W$, (i) and (ii) follow from Proposition
\ref{prop0}-\ref{prop2}. (iii), (iv) and (v) follow from (\ref{e12}), Observation \ref{ob}, and Proposition \ref{prop0}-\ref{prop2}.
\end{proof}

\begin{lemma}\label{typeAU}
 For any $u^*\in V(G)$, we have
 \vspace{-9pt}
\begin{enumerate}[(\romannumeral1)]
\setlength{\itemsep}{-1ex}
\item $u^*$ is incident with at most one nonremovable edge of type I;
\item $u^*$ is incident with at most two nonremovable edges of type II.
\end{enumerate}
\end{lemma}
\begin{proof}
Without loss of generality, we may assume that $u^*\in U$.

(i) Suppose to the contrary that $u^*$ is incident with two nonremovable edges of type I, say $u^*w_1$ and $u^*w_2$. 
By Lemma \ref{leme1e2e}, there is a barrier $S_i$ of $G-u^*w_i$ that contains $V(e_1)$, $i=1,2$. Thus $|S_1\cap S_2|\geq2$, contradicting  Lemma \ref{lem-barrier-cap}. (i) holds.

(ii) Suppose to the contrary that $u^*$ is incident with three nonremovable edges of type II, say $u^*w_1,u^*w_2$ and $u^*w_3$. Then  $u^*w_i$ ($1\leq i\leq 3$) is nonremovable in both $G$ and  $H$. By Lemma \ref{MC-nonremovable}, there exist partitions
$(U_0',U_1',U_2',U_3')$ of $U$ and $(W_0',W_1',W_2',W_3'$) of $W$, such that $u^*\in U_0'$, and for $i\in \{1,2,3\}$:
(a) $|U_i'|=|W_i'|$, (b) $w_i\in W_i'$, and (c) $N_H(W_i')=U_i'\cup \{u^*\}$; in particular, $u^*w_i$ is the only edge between $U_0'$ and $W_i'$.
Note that at least one of $W_1',W_2'$ and $W_3'$ contain no ends of $e_2$, say $W_3'$. Then $N(W_3')=U_3'\cup \{u^*\}$. By Lemma \ref{lem-barrier2},  $U_3'\cup \{u^*\}$ is a nontrivial barrier of $G$, a contradiction. (ii) holds.
\end{proof}

\section{Proof of Theorem \ref{th:main1}}

We first present a lemma, the proof of which will be given later.
\begin{lemma}\label{lem-2-edge}Let $u\in U\backslash V(e_1)$, and
$uy$ and  $uw$  be two nonremovable edges in $G$ such that $uy$ is of type I and $uw$ is of type II. Let $w'$ be a neighbour of $u$ other than $y$ and $w$. Then either $uw'$ is removable in $G$, or $uw'$ is of type II, $e_2= ww'$ and $d(u)=d(w)=d(w')=3$.
\end{lemma}

{\em{Proof of Theorem \ref{th:main1}.}}
Assume $u$ is a vertex in $U\backslash V(e_1)$  that is   incident with more than two nonremovable edges of $G$. Then, by Lemma \ref{typeAU}, $u$  is incident with exactly one edge of type I and two edges of type II.
We adopt the notational conventions stated
in Lemma \ref{lem-2-edge}.
Then  $uy$ is of type I, and $uw$ and $uw'$ are of type II.
By Lemma \ref{lem-2-edge}, we have $e_2=ww'$ and $d(u)=d(w)=d(w')=3$. Therefore, $u$ lies in the triangle containing $e_2$ and each vertex in this triangle is of degree three.

If there exist two vertices $u_1$ and $u_2$ in $U\backslash V(e_1)$ such that $u_1$ and $u_2$ are incident with more than two nonremovable edges of $G$, then we get two triangles $T_1$ and $T_2$ that contain  $u_1$ and $u_2$, respectively, and have the edge $e_2$ in common.
By Lemma \ref{lem-4-cycle}, for  $i=1, 2$,  $u_i$ is incident with a removable edge of $H$ in $T_i$ that is of type I in $G$, contradicting the fact that it is of type II.
Thus  $U\backslash V(e_1)$ contains  at most one vertex that is incident with more than two nonremovable edges of $G$.

The same result is also true for vertices in $W\backslash V(e_2)$.  So the theorem follows.
$\hfill\square$

\vspace{2mm}

Now we turn to the proof of Lemma \ref{lem-2-edge}.
By Corollary \ref{cor-P} (i) and (ii), $u\in U_2\cup U_4$, $w\in W_1\cup W_3$ and $\{uw\}=E[U_2\cup U_4, W_1\cup W_3]$. Thus $\{y,w'\}\subset W_2\cup W_4$.
Since $uy$ is of type I, Lemma \ref{leme1e2e} implies that there exists partitions $(A_1,A_2)$ of $U$ and $(B_1,B_2)$ of $W$ such
that $|B_1|=|A_1|+1$, $|A_2|=|B_2|+1$,  $\{e_1\}=E(A_2)$, $\{e_2\}=E(B_1)$ and $\{uy\}=E[A_1,B_2]$. Furthermore, $B_1$ and $A_2$ are barriers of $G-uy$.
Note that
\begin{equation}
u\in A_1\cap(U_2\cup U_4),  w\in B_1\cap(W_1\cup W_3), y\in B_2\cap (W_2\cup W_4),  w'\in B_1\cap (W_2\cup W_4),  \label{uyww'}
\end{equation}
\begin{equation}
N(A_1\setminus  \{u\})\subseteq B_1 \mbox{ and } N(B_2\setminus  \{y\})\subseteq A_2.   \label{AB}
\end{equation}
 Combining  (\ref{UW34}) and Corollary \ref{cor-P}(iii), we have
\begin{equation}N((U_i\cap A_1)\backslash\{u\})\subseteq(W_i\cup W_2)\cap B_1, i=3,4.\label{U4A1}\end{equation}

\subsection{Properties}
\begin{lemma}\label{lem-w'U''}
If $E[w',U_1\cup U_3]\neq \emptyset$, then $uw'$ is removable in $G$.
\end{lemma}

\begin{proof}
Suppose to the contrary that $uw'$ is nonremovable in $G$.
Since $uy$ is of type I, Lemma \ref{typeAU}(i) implies that $uw'$ is of type II, so $uw'$ is nonremovable in $H$.
Then, there exists an edge $f$ of $H-uw'$ such that each perfect matching of $H$ containing $f$ contains $uw'$. Let $M$ be a perfect matching  of $H$ containing both $f$ and $uw'$.
Let $M_1$ be a perfect matching of $H$  containing $uy$.
Then $uw$ is neither in $M$ nor in $M_1$.
Recall that  $\{uw\}=E[U_2\cup U_4,W_1\cup W_3]$ and $\{y,w'\}\subseteq W_2\cup W_4$.
By Corollary \ref{cor-P}(v), $|U_2\cup U_4|=|W_2\cup W_4|$. Let $X=U_2\cup U_4\cup W_2\cup W_4$. Then $M\cap \partial(X)=\emptyset$ and $M_1\cap \partial(X)=\emptyset$.
Since $f\in M$,  $f\in E(X)$ or $f\in E(\overline{X})$.
If $f\in E(\overline{X})$, then $(M\cap E(\overline{X}))\cup (M_1\cap E(X))$ is a perfect matching of $H$ that contains $f$ but does not  contain $uw'$, a contradiction.
So $f\in E(X)$. Recall that $E[w',U_1\cup U_3]\neq \emptyset$.
Let $z'$ be a neighbour of $w'$ in $U_1\cup U_3$ and $M_2$  a perfect matching of $H$ that contains $w'z'$.
Then $uw\in M_2$, so $M_2\cap \partial(X)=\{w'z',uw\}$. It follows that
$((M\cap E(X))\backslash\{uw'\})\cup (M_2\cap (E(\overline{X})\cup \partial(X)))$ is a perfect matching of $H$ that contains $f$ but does not contain $uw'$, a contradiction. The assertion follows.
\end{proof}

\begin{prop}\label{P-UAWB}
(i) $|W_2\cap B_1|\leq1$, $|U_1\cap A_2|\leq1$ and $W_1\subseteq  B_1$; (ii) If $\{e_1\}=E[U_1,U_2]$, then $|U_1\cap A_2|=1$; (iii) If $U_2\neq \emptyset$, then $W_2\cap B_1=\{w'\}$, $U_2\cap A_1=\{u\}$ and $d(u)=3$.
\end{prop}

\begin{proof}
(i) Recall that $U_1\cup W_2$ is a barrier of $G-uw$, and  $B_1$ and $A_2$ are barriers of $G-uy$. By Lemma \ref{lem-barrier-cap},  $|W_2\cap B_1|\leq1$ and $|U_1\cap A_2|\leq1$.
By (\ref{uyww'}), $w\notin B_2$ and $y\notin W_1$.
Recall that $\{e_2\}=E(B_1)$.
By (\ref{e*}) and (\ref{AB}), $N(W_1\cap B_2)\subseteq U_1\cap A_2$. Since $G$ is 3-connected, $W_1\cap B_2=\emptyset$, so $W_1\subseteq  B_1$. Therefore, (i) holds.

(ii) The result follows from  (i) and the fact that $\{e_1\}=E(A_2)\cap E[U_1,U_2]$.

(iii) If $U_2\neq \emptyset$, by (\ref{uyww'}) and Corollary \ref{cor-P}(i), $u\in U_2\cap A_1$. By (\ref{e*}) and (\ref{uyww'}), the neighbour $w'$ of $u$ lies in $W_2\cap B_1$. As $|W_2\cap B_1|\leq1$, we have $W_2\cap B_1=\{w'\}$, so $N(u)=\{y,w,w'\}$, i.e., $d(u)=3$.
By (\ref{e*}) and (\ref{AB}), $N((U_2\cap A_1)\backslash\{u\})\subseteq W_2\cap B_1$. Since $G$ is 3-connected, $(U_2\cap A_1)\backslash\{u\})=\emptyset$, so $U_2\cap A_1=\{u\}$. Therefore,
(iii) holds.
\end{proof}

\begin{prop}\label{P-UAWB+}
Assume that  $U_2=\emptyset$ and $uw'$ is nonremovable in $G$. Then $w'\in W_4\cap B_1$, $|W_4\cap B_1|\geq|U_4\cap A_1|$,  and $|(W_2\cup W_3)\cap B_1|\geq|U_3\cap A_1|+2$ if $W_3\neq \emptyset$.
\end{prop}

\begin{proof}
By Lemma \ref{lem-w'U''}, $E[w',U_1\cup U_3]=\emptyset$. By (\ref{e*}), $E[W_2,U_1]\neq\emptyset$. By  Corollary \ref{cor-P}(i),  $|W_2|=1$, so  $w'\notin W_2$. By (\ref{uyww'}), $w'\in W_4\cap B_1$.
By  Corollary \ref{cor-P}(iv),  if $W_3\neq \emptyset$, then $\{e_2\}=E(W_3)$. Recall that $\{e_2\}=E(B_1)$. Then $|W_3\cap B_1|\geq 2$.
By (\ref{U4A1}) and   Corollary \ref{cor-P}(i), $N(U_3\cap A_1)\subseteq (W_3\cup W_2)\cap B_1$ and $N((U_4\cap A_1)\backslash\{u\})\subseteq (W_4\cup W_2)\cap B_1$.
By Lemma \ref{lem-near-N(S)}, $|(W_2\cup W_3)\cap B_1|\geq|U_3\cap A_1|+2$  if $W_3\neq \emptyset$, and   $|W_4\cap B_1|\geq|U_4\cap A_1|$ if $W_2\subseteq B_1$.
If $W_2\subseteq B_2$,  when $(U_4\cap A_1)\backslash\{u\}\neq\emptyset$, Lemma \ref{lem-MC-N(S)} implies that $|W_4\cap B_1|\geq|U_4\cap A_1|$, which is also true when $U_4\cap A_1=\{u\}$. The assertion follows.
\end{proof}

By Proposition \ref{P-UAWB}(i),  we proceed to consider the following two cases.

\subsection{The case $|U_1\cap A_2|=1$}
\begin{lemma}\label{lem1-U2e2}
Assume that $|U_1\cap A_2|=1$.
\vspace{-9pt}
\begin{enumerate}[(\romannumeral1)]
\setlength{\itemsep}{-1ex}
\item If $U_2=\emptyset$ or $W_3\neq\emptyset$, then  $uw'$ is removable in $G$.
\item If $U_2\neq\emptyset$  and $W_3=\emptyset$, then $uw'$ is removable in $G$, or $uw'$ is of type II, $e_2= ww'$ and $d(u)=d(w)=3$.
\end{enumerate}
\end{lemma}

\begin{proof}
Since $|U_1\cap A_2|=1$, $|U_1\cap A_1|=|U_1|-1$. By Proposition \ref{P-UAWB}(i),   $W_1=W_1\cap B_1$. By Proposition \ref{prop0}-\ref{prop2}, $|W_1|\leq|U_1|\leq |W_1|+1$, so $|W_1|-1\leq|U_1\cap A_1|\leq |W_1\cap B_1|$. In particular, if  $|W_1|=|U_1|$, then $|W_1\cap B_1|=|U_1\cap A_1|+1$. We now show the following claim.

{\bf Claim.} If $U_2\neq \emptyset$ and $uw'$ is nonremovable in $G$, then $|W_2\cap B_1|=|U_2\cap A_1|=1$, $|W_4\cap B_1|\geq |U_4\cap A_1|$, and $|W_3\cap B_1|\geq|U_3\cap A_1|+2$ if $W_3\neq \emptyset$.

Since $U_2\neq \emptyset$, by Proposition \ref{P-UAWB}(iii), $U_2\cap A_1=\{u\}$ and $W_2\cap B_1=\{w'\}$. Then $u\notin U_3\cup U_4$.
By Lemma \ref{lem-w'U''}, $E[W_2\cap B_1,U_1\cup U_3]=\emptyset$. By (\ref{U4A1}), $N(U_3\cap A_1)\subseteq W_3\cap B_1$ and $N(U_4\cap A_1)\subseteq (W_4\cap B_1)\cup \{w'\}$.
By Lemma \ref{lem-MC-N(S)}, $|W_4\cap B_1|\geq |U_4\cap A_1|$.
If $W_3\neq\emptyset$, by Corollary \ref{cor-P}(iv), $\{e_2\}=E(W_3)$. Recall that $\{e_2\}=E(B_1)$. Then $|W_3\cap B_1|\geq 2$.
By Lemma \ref{lem-near-N(S)}, $|W_3\cap B_1|\geq|U_3\cap A_1|+2$.  Claim holds.

(i) Suppose to the contrary that $uw'$ is nonremovable  in $G$.
We will show   $|B_1|\geq |A_1|+2$, contradicting the fact that $|B_1|=|A_1|+1$.
Recall that $|W_1\cap B_1|\geq |U_1\cap A_1|$.
If $U_2\neq\emptyset$ and  $W_3\neq\emptyset$, by the above claim, we have $|B_1|=\sum_{i=1}^4|W_i\cap B_1|\geq\sum_{i=1}^4|U_i\cap A_1|+2=|A_1|+2$.
If $U_2=\emptyset$, by Proposition \ref{P-UAWB+}, $|W_4\cap B_1|\geq|U_4\cap A_1|$, and $|(W_2\cup W_3)\cap B_1|\geq|U_3\cap A_1|+2$ if $W_3\neq \emptyset$.
Then $|B_1|\geq|A_1|+2$.
If $W_3= \emptyset$, by Corollary \ref{cor-P}(iv),  $U_3=\emptyset$, $\{e_2\}=E[W_1,W_2]$ and $|U_1|=|W_1|$. Then $|W_1\cap B_1|=|U_1\cap A_1|+1$.
Recall that $\{e_2\}=E(B_1)$. By Corollary \ref{cor-P}(i), $|W_2|=1$, so $W_2\subseteq B_1$, i.e., $|W_2\cap B_1|=1$.
Consequently, $|B_1|\geq |A_1|+2$.
(i) holds.

(ii) Since $U_2\neq\emptyset$  and $W_3=\emptyset$, by Corollary \ref{cor-P},  $u\in U_2$, $\{e_2\}=E[W_1, W_2]$, $|U_1|=|W_1|$, $U_3=\emptyset$ and $w\in W_1$.  
By Proposition \ref{P-UAWB}(iii), $d(u)=3$ and $W_2\cap B_1=\{w'\}$, so $w'\in V(e_2)$.  Assume that   $uw'$ is nonremovable in $G$. By Lemma \ref{typeAU}, $uw'$ is of type II.
Suppose to the contrary that  $e_2\neq ww'$. Then $W_1$ contains at least two vertices, $w$ and one end of $e_2$, i.e.,  $|W_1|\geq2$.
Since $|U_1|=|W_1|$, we have $|W_1\cap B_1|=|U_1\cap A_1|+1$.
Since $|B_1|=|A_1|+1$, by the above claim, we have $|W_4\cap B_1|=|U_4\cap A_1|$.
Recall that $\{e_1\}=E(A_2)$.
Since $u\notin U_4$, by (\ref{U4A1}), $N(U_4\cap A_1)\subseteq (W_4\cap B_1)\cup \{w'\}$. By Lemma \ref{lem-near-N(S)}, $W_4\cap B_1=\emptyset$.
By Lemma \ref{lem-w'U''}, $E[W_2\cap B_1,U_1\cup U_3]=\emptyset$, so  $N(U_1\cap A_1)\subseteq W_1\cap B_1= W_1$.
By Lemma \ref{lem-barrier2},  $W_1$ is a nontrivial barrier of $G$.  This  contradiction implies that $e_2=ww'$. By (\ref{e*}), $N(W_1\backslash\{w\})\subseteq U_1$. Since $|U_1|=|W_1\backslash\{w\}|+1$, Lemma \ref{lem-barrier2} implies that $U_1$ is a barrier of $G$. By Lemma \ref{lem-barrier1}(i), $|U_1|=1$, so $W_1=\{w\}$. Then $w$ has exactly three neighbours $u$, $w'$ and the vertex in $U_1$, i.e., $d(w)=3$. (ii) holds.
\end{proof}

\subsection{The case $U_1\cap A_2=\emptyset$}

If   $U_1\cap A_2=\emptyset$, then $U_1\subseteq A_1$, i.e., $U_1\cap A_1=U_1$. Recall that $W_1\cap B_1=W_1$. By Corollary \ref{cor-P}(i),  $u\notin U_1$ and,  by (\ref{AB}), $N(U_1)\subseteq N(A_1\setminus\{u\})\subseteq B_1$.

\begin{lemma}\label{lem0-U_2=2}
If  $U_1\cap A_2=\emptyset$ and $U_2\neq\emptyset$, then  $uw'$ is removable in $G$.
\end{lemma}

\begin{proof}
Suppose to the contrary that $uw'$ is nonremovable in $G$.
Since $U_2\neq\emptyset$, by Proposition \ref{P-UAWB}(iii),  $W_2\cap B_1=\{w'\}$.  By (\ref{e*}), $E[U_1,W_2]\neq\emptyset$. Since $N(U_1)\subseteq B_1$, we have $E[U_1,W_2\cap B_1]\neq\emptyset$, i.e., $E[U_1,w']\neq\emptyset$, contradicting Lemma \ref{lem-w'U''}.
\end{proof}

To deal with the case $U_1\cap A_2=\emptyset$ and $U_2=\emptyset$, we need to consider the properties with respect to $uw'$.
If $uw'$ is of type II, assume that $(U_1',U_2', U_3',U_4')$ and  $(W_1',W_2',W_3',W_4')$  are  partitions of $U$ and $W$, respectively, with respect to $uw'$ defined as those with respect to $uw$ in Section 2.
Then all the properties with respect to $uw$ are also true with respect to $uw'$.
If $U_2=\emptyset$, by Proposition \ref{prop0}(i),  $\omega\geq 1$.
We say the edge $uw$ is of {\em type 1} if  $\omega=1$, and of {\em type 2} if $\omega=2$. Analogously, we may define the type of the edge $uw'$ when $U_1'\cap A_2=\emptyset$ and $U_2'=\emptyset$.

\begin{lemma}\label{lem0-U2=0}
Assume  that  $U_1\cap A_2=\emptyset$, $U_2=\emptyset$ and $uw'$ is nonremovable in $G$. Then $uw'$ is  neither  of type 1 nor of type 2.
\end{lemma}

The proof of Lemma \ref{lem0-U_2=1} requires Lemma \ref{lem0-U2=0}, whose proof is presented in the next subsection.

\begin{lemma}\label{lem0-U_2=1}
If  $U_1\cap A_2=\emptyset$ and $U_2=\emptyset$, then  $uw'$ is removable in $G$.
\end{lemma}

\begin{proof}
Suppose to the contrary that $uw'$ is nonremovable in $G$.
By Proposition \ref{P-UAWB+}, $w'\in W_4$.
Recall that $N(u)=\{w,w',y\}$,  $uy$ is of type I and $uw$ is of type II.  By Lemma \ref{typeAU}(i), $uw'$ is of type II.
So $|U_1'\cap A_2|\leq 1$. By Lemma \ref{lem0-U2=0}, $uw'$ is  neither  of type 1 nor of type 2. Then  $|U_1'\cap A_2|=1$ or $U_2'\neq\emptyset$.
By Corollary \ref{cor-P}(iv), $\{e_2\}=E(W_3)\cup E[W_1,W_2]$.
Since $w'\in W_4$, $w'\notin V(e_2)$, so $e_2\neq ww'$. If $|U_1'\cap A_2|=1$, Lemma \ref{lem1-U2e2} implies that $uw$ is removable in $G$, a contradiction. If $U_1'\cap A_2=\emptyset$ and $U_2'\neq\emptyset$, Lemma \ref{lem0-U_2=2} implies that $uw$ is removable in $G$, also a contradiction.
\end{proof}

{\em{Proof of Lemma \ref{lem-2-edge}.}}
By Lemmas \ref{lem1-U2e2}, \ref{lem0-U_2=2} and \ref{lem0-U_2=1}, we only need to show $d(w')=3$ when $uw'$ is of type II. 
 In this case, noticing $uw$ is of type II, with $w'$ playing the role of $w$, we have $d(w')=3$ by
 Lemmas \ref{lem1-U2e2}, \ref{lem0-U_2=2} and \ref{lem0-U_2=1} again.
$\hfill\square$

\subsection{Proof of Lemma \ref{lem0-U2=0}}

In  this subsection, we  show Lemma \ref{lem0-U2=0}  by contradiction. Suppose that $uw'$ is of type 1 or type 2.
Then $uw'$ is of type II,  $U_1'\cap A_2=\emptyset$, $U_2'=\emptyset$ and $|W_2'|=1$.
We may let $W_2'=\{w_2'\}$.
We also have (a) $\{e_1\}=E(U_4')$, (b) if $W_1'\neq \emptyset$, then $\{uw'\}=E[U_4',W_1']$; otherwise, $\{uw'\}=E[U_4',W_3']$, and (c)  if $W_3'\neq \emptyset$, then $\{e_2\}=E(W_3')$; otherwise,  $\{e_2\}=E[W_1',W_2']$.
Recall that $uw$ is of type II. Since $U_1\cap A_2=\emptyset$ and $U_2=\emptyset$, $uw$ is also of type 1 or type 2, so $|W_2|=1$. Let $W_2=\{w_2\}$. 
Note that $u\in U_4\cap U_4'\cap A_1$ and $\{e_1\}=E(U_4)\cap E(U_4')$. 
The following properties with respect to $uw$ are also true with respect to $uw'$.

\begin{prop}\label{P-1122}
$|U_4\cap A_1|=|W_4\cap B_1|$,  $A_2\subseteq U_4$ and $B_2\subseteq W_4$.
\end{prop}

\begin{proof} Recall that $W_1\subseteq B_1$, $U_1\subseteq A_1$, $U_2=\emptyset$ and $N(U_1)\subseteq B_1$. By (\ref{e*}), $E[U_1, w_2]\neq\emptyset$, so $w_2\in B_1$, i.e., $|W_2\cap B_1|=1$.
By Proposition \ref{P-UAWB+}, $|W_4\cap B_1|\geq|U_4\cap A_1|$, and $|(W_2\cup W_3)\cap B_1|\geq|U_3\cap A_1|+2$ if $W_3\neq \emptyset$.
If $uw$ is of type 1, by Proposition \ref{prop1}(ii), $|W_1|=|U_1|$, so $|W_1\cap B_1|=|U_1\cap A_1|$. By Observation \ref{ob}, $U_3=W_3=\emptyset$, so $A_2\subseteq U_4$ and $B_2\subseteq W_4$.
Since $|B_1|=|A_1|+1$, we have $|W_4\cap B_1|=|U_4\cap A_1|$.
If $uw$ is of type 2, by (\ref{e12}) and Proposition \ref{prop2}, $W_3\neq \emptyset$ and $|W_1|=|U_1|-1$, so $|W_1\cap B_1|=|U_1\cap A_1|-1$.
Since $|B_1|=|A_1|+1$, we have $|W_4\cap B_1|=|U_4\cap A_1|$ and  $|W_3\cap B_1|=|U_3\cap A_1|+1$. Note that $|W_3|=|U_3|+1$. Then $|W_3\cap B_2|=|U_3\cap A_2|$.
Recall that $w\in B_1$, $y\in B_2\cap (W_2\cup W_4)$ and $N(U_1)\subseteq B_1$. By (\ref{UW34}) and (\ref{AB}), $N(W_3\cap B_2)\subseteq U_3\cap A_2$.
By Lemma \ref{lem-MC-N(S)}, $W_3\cap B_2=\emptyset$, so $U_3\cap A_2=\emptyset$, i.e., $ W_3\subseteq B_1$ and $U_3\subseteq A_1$. Thus $A_2\subseteq U_4$ and $B_2\subseteq W_4$.
\end{proof}

Recall that $\{w_2\}\cup U_1$ and $\{w_2'\}\cup U_1'$ are barriers of $G-uw$ and $G-uw'$, respectively. Lemma \ref{lem-barrier-cap} implies that
\begin{equation}\label{barriers-uww'} |(\{w_2\}\cup U_1)\cap(\{w_2'\}\cup U_1')|\leq 1. \end{equation}
Since $u\in U_4$ and $\{e_1\}=E(U_4)$, by (\ref{UW34}),
\begin{equation}\label{U3Uj'} N(U_3\cap U_j')\subseteq (W_3\cup \{w_2\})\cap(W_j'\cup \{w_2'\}), \mbox{ where } j=3,4. \end{equation}

\begin{prop}\label{P-UAWB+3}
(i)  $W_1\cap W_1'=\emptyset$; and
(ii) if $w_2\in W_1'\cup W_3'$ and $w_2'\in W_1\cup W_3$, 
then $|W_4\cap W_4'\cap B_1|\geq|U_4\cap U_4'\cap A_1|+1$.
\end{prop}

\begin{proof} By Proposition \ref{P-UAWB+},  $w'\in W_4$ and $w\in W_4'$.

(i) By (\ref{barriers-uww'}), $|U_1\cap U_1'|\leq 1$.
Note that  $w,w'\notin W_1\cap W_1'$. By (\ref{e*}), $N_{G-e_2}(W_1\cap W_1')\subseteq U_1\cap U_1'$.
As $G$ is 3-connected,  (i) follows.

(ii) Note that $w_2\neq w'$ and $w_2'\neq w$.
Since $w_2\in W_1'\cup W_3'$ and $w_2'\in W_1\cup W_3$, by (\ref{e*}) and (\ref{UW34}),  $w_2$ and $w_2'$ have no neighbours in $U_4\cap U_4'$. Recall that $u\in U_4\cap U_4'\cap A_1$. By Proposition \ref{P-1122}, $A_2\subseteq (U_4\cap U_4')\setminus \{u\}$. By (\ref{UW34}), $N(A_2)\subseteq N((U_4\cap U_4')\backslash\{u\})\subseteq  W_4\cap W_4'$.
By (\ref{AB}), $N(B_2\backslash\{y\})\subseteq A_2$.
Then $N(A_2\cup (B_2\backslash\{y\}))\subseteq
\{y\}\cup (W_4\cap W_4'\cap B_1)$. Since $G$ is 3-connected, $|W_4\cap W_4'\cap B_1|\geq2$.
By (\ref{UW34}) and (\ref{AB}), $N((U_4\cap U_4'\cap A_1)\backslash\{u\})\subseteq W_4\cap W_4'\cap B_1$, (ii) follows from Lemma \ref{lem-near-N(S)}.
\end{proof}

\begin{prop}\label{lem-uwuw'-1}
Neither  $uw$ nor $uw'$ is of type 2.
\end{prop}

To complete the proof of Lemma \ref{lem0-U2=0}, we need Proposition \ref{lem-uwuw'-1}, the proof of that is a main burden.
We present it later. From Proposition \ref{lem-uwuw'-1}, we see that both $uw$ and $uw'$ are of type 1.
By Proposition \ref{prop1} and Observation \ref{ob}, $\{e_2\}=E[W_1,W_2]\cap E[W_1',W_2']$ and $U_3=W_3=\emptyset$.
Then each of $w_2$ and $w_2'$ is an end of $e_2$.
By Proposition \ref{P-UAWB+3}(i), $W_1\cap W_1'=\emptyset$.
If $w_2=w_2'$, then the other end of $e_2$ lies in $W_1\cap W_1'$, a contradiction. 
Thus $w_2\neq w_2'$. Then  $e_2=w_2w_2'$, so $w_2'\in W_1$ and $w_2\in W_1'$.
The former implies that  $|W_1\cap W_2'|=1$. Proposition \ref{P-UAWB+3}(ii) implies that $|U_4\cap U_4'\cap A_1|\leq |W_4\cap W_4'\cap B_1|-1$.
By Proposition \ref{P-1122}, $|U_4'\cap A_1|=|W_4'\cap B_1|$, $U_1\subseteq A_1$ and $W_1\subseteq B_1$. Note that $(U_1, U_4)$ is a partition of $U$ and  $(W_1, W_2, W_4)$ is a partition of $W$. Since $w_2\in W_1'$,
we have $|U_1\cap U_4'\cap A_1|\geq|W_1\cap W_4'\cap B_1|+1$, so $|U_1\cap U_4'|\geq|W_1\cap W_4'|+1$. By Proposition \ref{P-UAWB+}, $w'\in W_4$, so $w_2\neq w'$. By (\ref{e*}), $N(w_2)\subseteq U_1'\cup\{w_2'\}$.
Since $E[w_2, U_1]\neq\emptyset$ and  $d(w_2)\geq 3$, we have $U_1\cap U_1'\neq\emptyset$, so $|U_1\cap U_1'|=1$.
Consequently, $|U_1|\geq |W_1|+1$,
contradicting the fact that $|U_1|=|W_1|$.
So Lemma \ref{lem0-U2=0} holds.

\subsection{Proof of Proposition \ref{lem-uwuw'-1}}
Assume, without loss of generality, that $uw$ is of type 2. We first present some basic properties.


\begin{prop}\label{prop1'} Assume that  $|U_1|\geq2$. Then  $|U_1\cap U_1'|+|U_1\cap U_4'|=|W_1\cap W_4'|+1$ and, if $w_2'\notin W_1$, then $U_1\cap U_3'=W_1\cap W_3'=\emptyset$; otherwise,  $|U_1\cap U_1'|=1$, $|U_1\cap U_3'|=|W_1\cap W_3'|+1$  and $w_2\in W_1'\cup W_3'$.
\end{prop}

\begin{proof}
Since $uw$ is of type 2, $|U_1|=|W_1|+1$, so $W_1\neq \emptyset$. By Corollary \ref{cor-P}(ii) and Proposition \ref{P-UAWB+}, $w\in W_1\cap W_4'$. Note that $\{e_2\}=E(W_3)$.
By (\ref{e*}) and (\ref{UW34}), $N(W_1\cap W_4')\subseteq (U_1\cap (U_4'\cup U_1'))\cup\{u\}$ and $N(W_1\cap W_3')\subseteq U_1\cap (U_1'\cup U_3')$. For the former case, since $d(w)\geq3$, we have $|U_1\cap (U_4'\cup U_1')|\geq2$. By Lemma \ref{lem-near-N(S)}, $|U_1\cap (U_4'\cup U_1')|\geq |W_1\cap W_4'|+1$.
For the latter, if $W_1\cap W_3'\neq\emptyset$, Lemma \ref{lem-MC-N(S)} implies that $|U_1\cap (U_1'\cup U_3')|\geq|W_1\cap W_3'|+1$. By (\ref{barriers-uww'}), $|U_1\cap U_1'|\leq 1$, so $|U_1\cap U_3'|\geq|W_1\cap W_3'|$, which is also true when $W_1\cap W_3'=\emptyset$.
Recall that $W_1\cap W_1'= \emptyset = U_2=U_2'$.

Assume that $w_2'\notin W_1$. Then $W_1\cap W_2'=\emptyset$.
Since $|U_1|=|W_1|+1$, we have $\sum_{i=1}^4|U_1\cap U_i'|=\sum_{i=1}^4|W_1\cap W_i'|+1$.
Then $|U_1\cap (U_4'\cup U_1')|=|W_1\cap W_4'|+1$ and $|U_1\cap U_3'|=|W_1\cap W_3'|$, the latter implies that  $|U_1\cap (U_1'\cup U_3')|\leq|W_1\cap W_3'|+1$.
Since $W_1\cap W_3'$ is an independent set, Lemma \ref{lem-barrier2} implies that  $U_1\cap (U_1'\cup U_3')$ is a barrier of $G$. By Lemma \ref{lem-barrier1}(i), $U_1\cap (U_1'\cup U_3')$ is a singleton. As $G$ is 3-connected, $W_1\cap W_3'=\emptyset=U_1\cap U_3'$.

Assume that $w_2'\in W_1$. Then $uw'$ is of type 2, otherwise, $w_2'\in V(e_2)$, so $w_2'\in W_3$, a contradiction. So $\{e_2\}=E(W_3')$ and $|W_3'|=|U_3'|+1$.
If $E[w_2',U_3']=\emptyset$, then $N(U_3')\subseteq W_3'$. By Lemma \ref{lem-barrier2}, $W_3'$ is a nontrivial barrier of $G$, a contradiction. Therefore, $E[w_2',U_3']\neq\emptyset$.
Recall that $E[w_2',U_1']\neq\emptyset$ and $w_2'\neq w$.
Since $w_2'\in W_1$ and $N(W_1\backslash\{w\})\subseteq U_1$, we have $N(w_2')\subseteq U_1$.
So $U_1\cap U_1'$ and $U_1\cap U_3'$ are nonempty, the former implies that $|U_1\cap U_1'|=1$.
By Lemma \ref{lem-near-N(S)},
$|U_1\cap (U_1'\cup U_3')|\geq|W_1\cap W_3'|+2$, so $|U_1\cap U_3'|\geq|W_1\cap W_3'|+1$.
Since $|U_1|=|W_1|+1$ and  $|W_1\cap W_2'|=1$, by the same reason as the case $w_2'\notin W_1$, we have $|U_1\cap (U_4'\cup U_1')|=|W_1\cap W_4'|+1$ and $|U_1\cap U_3'|=|W_1\cap W_3'|+1$.

We now show that $w_2\in W_1'\cup W_3'$. If not,  $w_2\in \{w_2'\}\cup W_4'$.  By (\ref{barriers-uww'}),  $w_2\neq w_2'$, so $w_2\in W_4'$.
If $|U_1'|=1$,  then $W_1'=\emptyset$. Combining $|U_1\cap U_1'|=1$, we have $W_3\cap W_1'=U_3\cap U_1'=\emptyset$, which is also true when  $|U_1'|\geq2$ (as above) since $w_2\notin W_1'$.
Recall that  $w_2'\in W_1$. By (\ref{U3Uj'}),  $N(U_3\cap U_3')\subseteq W_3\cap W_3'$ and $N(U_3\cap U_4')\subseteq\{w_2\}\cup(W_3\cap W_4')$.
Because $\{e_2\}=E(W_3)\cap E(W_3')$, Lemma \ref{lem-near-N(S)} implies that $|W_3\cap W_3'|\geq|U_3\cap U_3'|+2$.
If $U_3\cap U_4'\neq\emptyset$, Lemma \ref{lem-MC-N(S)} implies that $|W_3\cap W_4'|\geq|U_3\cap U_4'|$, which is also true when $U_3\cap U_4'=\emptyset$.
Note that $W_3\cap W_2'=\emptyset$.
We have $|W_3|\geq|U_3|+2$, contradicting the fact that $|W_3|=|U_3|+1$.
\end{proof}

\begin{prop}\label{prop2'} Assume that $w_2'\in W_1\cup W_3$.
If $|U_1|\geq2$ and $|U_1\cap U_1'|=1$, or $|U_1|=1$ and $U_1\cap U_4'=\emptyset$, then $w_2\in W_4'$.
\end{prop}

\begin{proof}
Suppose to the contrary that $w_2\notin W_4'$. Then $W_2\cap W_4'=\emptyset$.
Since $w_2'\in W_1\cup W_3$, we have $w_2\neq w_2'$, so $w_2\in W_1'\cup W_3'$.
By Proposition \ref{P-UAWB+3}(ii),  $|W_4\cap W_4'\cap B_1|\geq|U_4\cap U_4'\cap A_1|+1$.
By (\ref{U3Uj'}),  $N(U_3\cap U_4')\subseteq W_3\cap(W_4'\cup \{w_2'\})$.
Lemma \ref{lem-MC-N(S)} implies that $|W_3\cap W_4'|\geq|U_3\cap U_4'|$.
We  assert that $|W_1\cap W_4'|=|U_1\cap U_4'|$, which  follows from Proposition \ref{prop1'} if $|U_1|\geq2$ and $|U_1\cap U_1'|=1$, and from the fact that $U_1\cap U_4'=\emptyset$ and $W_1=\emptyset$ if $|U_1|=1$.
By Proposition \ref{P-1122}, $U_1,U_3\subseteq A_1$ and $W_1,W_3\subseteq B_1$. Then $|W_4'\cap B_1|=\sum_{i=1}^4|W_4'\cap B_1\cap W_i|\geq \sum_{i=1}^4|U_4'\cap A_1\cap U_i|+1= |U_4'\cap A_1|+1$, contradicting the fact that $|W_4'\cap B_1|=|U_4'\cap A_1|$.
\end{proof}

If $uw'$ is of type 1, then  $w'\in W_1'$,  $\{e_2\}= E[w_2',W_1']$, and let $k=1$; otherwise, $w'\in W_1'$ or  $w'\in W_3'$, $\{e_2\}=E(W_3')$,  and let $k=3$.
Since $uw$ is of type 2, we have $\{e_2\}=E(W_3)$.
Thus, $W_3\cap W_k'$ contains at least one end of $e_2$, i.e., $W_3\cap W_k'\neq \emptyset$.
By Proposition \ref{P-UAWB+},  $w\in W_4'$ and $w'\in W_4$, so $w'\notin V(e_2)$ and $N_H(W_3\cap W_k')\subseteq (U_1\cup U_3)\cap(U_1'\cup U_k')$.
Lemma \ref{lem-MC-N(S)} implies that
\begin{equation} \label{U13W3}
|(U_1\cup U_3)\cap(U_1'\cup U_k')|\geq|W_3\cap W_k'|+1.
\end{equation}
By Lemma \ref{lem-w'U''}, $E[w',U_1\cup U_3]=\emptyset$, so $N(w')\subseteq (U_4\cap (U_1'\cup U_k'))\cup\{u\}$.
Since $d(w')\geq3$,  $|U_4\cap (U_1'\cup U_k')|\geq2$. Since
$N((W_4\cap W_k')\backslash\{w'\})\subseteq (U_1\cup U_4)\cap(U_1'\cup U_k')$, Lemma \ref{lem-near-N(S)} implies that
\begin{equation} \label{U14W4}
|(U_1\cup U_4)\cap(U_1'\cup U_k')|\geq|W_4\cap W_k'|+1.
\end{equation}

If $uw'$ is of type 1, then $k=1$ and $|U_1'|=|W_1'|$. By (\ref{U13W3}) and (\ref{U14W4}), $|W_1'|=|(U_1\cup U_4)\cap U_1'|+|(U_1\cup U_3)\cap U_1'|-|U_1\cap U_1'|\geq|W_4\cap W_1'|+|W_3\cap W_1'|-|U_1\cap U_1'|+2$.
Since $W_1\cap W_1'=\emptyset$, we have $|W_1'|=\sum_{i=2}^4|W_i\cap W_1'|$, so
$|W_2\cap W_1'|+|U_1\cap U_1'|\geq 2$.  Since $|U_1\cap U_1'|\leq 1$,
we have $|U_1\cap U_1'|=1$ and  $w_2\in W_1'$.
If $|U_1|=1$, then $U_1\cap U_4'=\emptyset$.
Since $w_2'\in V(e_2)$,  $w_2'\in W_3$.
By  Proposition \ref{prop2'},  $w_2\in W_4'$, a contradiction.
Then  $uw'$ is of type 2, so  $|W_3'|=|U_3'|+1$ and $k=3$.
Recall that $uw$ is  of type 2. Then $\{e_2\}=E(W_3)\cap E(W_3')$.
Let $Z_3=(U_3\cap U_3')\cup (W_3\cap W_3')$. Then  $|Z_3|\geq 2$.
Since $w,w'\notin W_3\cap W_3'$, $N(W_3\cap W_3')\subseteq(U_1\cup U_3)\cap(U_1'\cup U_3')$. Combining (\ref{U3Uj'}), we have
\begin{equation}\label{Z3}
N(Z_3)\subseteq(\{w_2\}\cap(\{w_2'\}\cup W_3'))\cup (\{w_2'\}\cap W_3)
\cup (U_1\cap(U_1'\cup U_3'))\cup (U_1'\cap U_3).
\end{equation}

If $|U_1|=1$ or $|U_1'|=1$, let $U_1=\{u_1\}$ or $U_1'=\{u_1'\}$, respectively.

\begin{prop}\label{prop3}
If $|U_1|\geq2$ and $w_2'\notin W_1$, or $|U_1|=1$, then $|U_1'|\geq2$. Analogously, if $|U_1'|\geq2$ and $w_2\notin W_1'$, or $|U_1'|=1$, then $|U_1|\geq2$.
\end{prop}

\begin{proof}
Suppose to the contrary that $|U_1'|=1$. Then $W_1'=\emptyset$ and $|(U_1\cup U_3\cup U_4)\cap U_1'|=1$.
First consider the case $|U_1|\geq2$ and $w_2'\notin W_1$.
By Proposition \ref{prop1'},  $U_1\cap U_3'=W_1\cap W_3'=\emptyset$. By (\ref{U13W3}) and (\ref{U14W4}), $|U_1\cap U_1'|+|U_4\cap U_3'|+|U_3\cap U_3'|\geq|W_4\cap W_3'|+|W_3\cap W_3'|+1$.
If $|U_1\cap U_1'|=0$, since  $|W_2\cap  W_3'|\leq 1$, we have $|U_3'|\geq|W_3'|$, a contradiction.
If $|U_1\cap U_1'|=1$, 
by (\ref{barriers-uww'}), $w_2'\neq w_2$, so $w_2'\in W_3\cup W_4$.
If $w_2'\in W_4$, by (\ref{Z3}) and the fact that $U_1\cap U_3'=\emptyset$, then $N(Z_3)\subseteq \{u_1', w_2\}$, contradicting the fact that $G$ is 3-connected. 
Thus $w_2'\in W_3$. By Proposition \ref{prop2'},  $w_2\in W_4'$, i.e., $|W_2\cap W_3'|=0$. Therefore, $|U_3'|\geq|W_3'|$, a contradiction.

Now assume that $|U_1|=1$. Then $W_1=\emptyset$.
If $w_2'=w_2$, since $G$ is 3-connected, by (\ref{Z3}), we have $N(Z_3)=\{w_2,u_1,u_1'\}$, $u_1\in U_3'$ and $u_1'\in U_3$. By (\ref{U13W3}) and (\ref{U14W4}), we have $|U_3\cap U_3'|\geq|W_3\cap W_3'|-1$ and $|(U_1\cup U_4)\cap U_3'|\geq|W_4\cap W_3'|+1$. Consequently,  $|U_3'|\geq|W_3'|$, a contradiction. So $w_2'\neq w_2$. Then $w_2'\in W_3\cup W_4$.
Recall that $w_2'u_1'\in E(G)$.
If $w_2'\in W_4$, then $u_1'\in U_4\cup U_1$.
By (\ref{Z3}), $N(Z_3)\subseteq\{u_1,w_2\}$, a contradiction. Then $w_2'\in W_3$, so $u_1'\in U_1\cup U_3$.
Analogously, $w_2\in W_3'$ and $u_1\in U_1'\cup U_3'$.
Then $U_1\cap U_4'=\emptyset$.
By Proposition \ref{prop2'}, $w_2\in W_4'$, a contradiction.
\end{proof}

We are ready to complete the proof of Proposition \ref{lem-uwuw'-1} by  distinguishing the following two cases to get contradictions. By Proposition \ref{P-UAWB+},  $w\in W_4'$, so $w_2'\neq w$.

{\bf Case 1.} $|U_1|\geq2$ and $|U_1'|\geq2$. By Proposition \ref{prop1'}, we  may claim  that  $U_1\cap U_3'=\emptyset$ if $w_2'\notin W_1$ and $U_1'\cap U_3=\emptyset$  if $w_2\notin W_1'$.
First suppose that $U_1\cap U_1'=\emptyset$. We assert that $w_2'\notin W_1$.  Otherwise,  $N(w_2')\subseteq N(W_1\backslash\{w\})\subseteq U_1$. Since $E[w_2', U_1']\neq\emptyset$,  $N(w_2')\cap (U_1\cap U_1')\neq \emptyset$, a contradiction.
Analogously, $w_2\notin W_1'$. By the above claim,  $U_1\cap U_3'=U_1'\cap U_3=\emptyset$.
By (\ref{Z3}),  $N(Z_3)\subseteq \{w_2,w_2'\}$,  contradicting the fact that $G$ is 3-connected.
Now suppose that $|U_1\cap U_1'|=1$.  By (\ref{barriers-uww'}), $w_2'\neq w_2$.
When $w_2'\in W_1\cup W_3$, Proposition \ref{prop2'} implies that  $w_2\in W_4'$, so $w_2\notin W_1'$.
If $w_2'\notin W_1$, by (\ref{Z3}),  $N(Z_3)\subseteq \{w_2'\}\cup (U_1\cap U_1')$,  a contradiction. Thus $w_2'\in W_1$.
By Proposition \ref{prop1'},   $w_2\in W_1'\cup W_3'$, a contradiction.
When $w_2'\in W_4$, analogously,  we deduce that $w_2\in W_1'$, so $w_2'\in W_1\cup W_3$, a contradiction.

{\bf Case 2.} At least one of $U_1$ and $U_1'$ is a singleton. Assume, without loss of generality,  that $|U_1'|=1$. By Proposition \ref{prop3},  $|U_1|\geq2$ and $w_2'\in W_1$.
Proposition \ref{prop1'} implies that $w_2\in W_1'\cup W_3'$. 
Recall that $w_2'\neq w$.
Since $N(w_2')\subseteq N(W_1\backslash\{w\})\subseteq U_1$ and $u_1'w_2'\in E(G)$, we have  $u_1'\in U_1$, so $|U_1'\cap U_1|=1$.
By Proposition \ref{prop2'}, $w_2\in W_4'$, a contradiction.

\section{Proof of Theorem \ref{th:main2}}

Let $G$ be a graph and $C=\partial(X)$ an edge cut of $G$.
We call $G/X$ and $G/\overline{X}$ the \emph{$C$-contractions} of $G$, and we say that $G$ is a \emph{splicing} of $G/X$ and $G/\overline{X}$ (at $x$ and $\overline x$), or $G$ is obtained by {\it splicing} $G/X$ and $G/\overline{X}$ (at $x$ and $\overline x$).
A graph $G$ is a \emph{tri-ladder} if it is $\overline{C_6}$ or it may be obtained from $\overline{C_6}$ by iterative splicings $K_4$ at vertices in triangles. 
More precisely, there exists a sequence of graphs $(G_0, G_1, \ldots, G_r)$ such that $G_0=\overline{C_6}$, $G=G_r$, and, for $1\leq i \leq r$,  $G_i$ is a splicing of $G_{i-1}$ and $K_4$ at vertices in triangles.
Note that each $G_i$ is a tri-ladder and $G_0$ $(=\overline{C_6})$  has exactly two disjoint triangles, say $T$ and $T_0$. For convenience, we may assume that the vertices in $T$ are never used to splice with $K_4$.
Let  $u,v$ and $w$ be three vertices of $T_0$. Let $w'$ be a vertex of $K_4$. Suppose that  $G_1$ is a splicing of $G_0$ and $K_4$ at $w$ and  $w'$.
Then $uv \in E(G_1)$, and $uv$ is referred to as a \emph{rung} of $G_i$ $(1\leq i \leq r)$. Moreover, $G_1$ has two disjoint triangles, one is $T$ and the other is $T_1=K_4-w'$.
Continuing in this way, we see that each $G_i$,  $1\leq i \leq r$, has one more rung than $G_{i-1}$, and has two disjoint triangles, one  is  $T$ and the other is denoted by $T_i$.
Then $G_i$ has $i$ rungs. For the rung $e$ in $G_i$ not in $G_{i-1}$, we refer to $i$ as the {\it rank} of $e$, denoted by $R(e)$.
Note that  the $r$ rungs form a matching.  Deleting them from $G$, the resulting graph has exactly three vertex-disjoint paths connecting vertices in $T$ and $T_r$, which are referred to as \emph{ridges} of $G$.
The four graphs depicted in Figure 1 are near-bipartite tri-ladders whose rungs are illustrated by bold lines.
\vspace{1mm}

Since $K_4$ and $\overline {C_6}$ are cubic bricks, by the following lemma,  a tri-ladder is a cubic brick.

\begin{lemma}[\cite{CLM05}]\label{lem:cubic-brick-splicing}
Any splicing of two cubic bricks is a cubic brick.
\end{lemma}

\begin{lemma}[\cite{CLM05}]\label{lem:mc-splicing}
Any splicing of two matching covered graphs is a matching covered graph.
\end{lemma}

The following lemma is an immediate consequence of Lemma \ref{lem:mc-splicing}, also see \cite{LFL22^{+}}.

\begin{lemma}\label{lem41}
Suppose that $G$ is a splicing of two  matching covered graphs $G_1$ and $G_2$ at $u_1$ and $u_2$, where $u_i\in V(G_i)$, $i=1,2$.  Then every removable edge in $G_i$ that is not incident with $u_i$ is removable in $G$.
\end{lemma}

%
%

For a bipartite graph with a perfect matching, by Hall's Theorem, we can obtain a characterization
of  nonadmissible edges of the graph. Using this characterization, the following lemma can be easily proved.

\begin{lemma}[\cite{LFL22^{+}}]\label{lem:near-bipartite}
Let $G$ be a 3-connected cubic nonbipartite graph. If $G$ has a pair of edges $e_1$ and $e_2$ such that $G-\{e_1,e_2\}$ is bipartite, then  $G$ is a near-bipartite graph with removable doubleton $\{e_1,e_2\}$.
\end{lemma}


\begin{lemma}[\cite{CLM99}]\label{lem:brace}
In a brace on six or more vertices, every edge is removable.
\end{lemma}

\begin{lemma}[\cite{KCLL20}]\label{lem:3-cut}
In a cubic matching covered graph, each 3-cut is a separating cut.
\end{lemma}

It is known that a cubic graph  is 3-connected if and only if it is 3-edge-connected.  A cubic graph is termed {\it essentially 4-edge-connected} if it is 2-edge-connected and  free of nontrivial 3-cuts. Note that every nontrivial 3-cut is a matching in a 3-connected graph. An edge cut $C$ of a matching covered $G$ is {\it good} if it is separating but not tight.

\begin{lemma}[\cite{KCLL20}]\label{lem:4-edge}
Every essentially 4-edge-connected cubic graph is either a brick or a brace.
\end{lemma}


\begin{lemma}\label{lem43}
 Let $G$ be a 3-connected cubic near-bipartite graph with removable doubleton $\{e_1,e_2\}$, and $H=G-\{e_1,e_2\}$ with a bipartition $(U,W)$ such that $V(e_1)\subseteq U$.
Suppose that $\partial(X)$ is a nontrivial 3-cut of $G$. Then $|X|$ is odd, and $G/\overline{X}$ and $G/X$ are 3-connected cubic graphs. Moreover, if $|X\cap U|\geq|X\cap W|$, then the following statements hold.
\vspace{-9pt}
\begin{enumerate}[(\romannumeral1)]
\setlength{\itemsep}{-1ex}
\item  $|X\cap U|=|X\cap W|+1$.
\item (\cite{NK19}) $\partial(X)$ is tight in $G$ if and only if $E[X\cap W, \overline{X}\cap U]=\emptyset$, one of $e_1$ and $e_2$ has both ends in $X\cap U$ or $\overline{X}\cap W$ (adjust notation so that $V(e_1)\subseteq X\cap U$), and $e_2$ has at least one end in $X\cap W$. In addition, $G/X$ is bipartite and $G/\overline{X}$ is a near-bipartite graph with removable doubleton $\{e_1,e_2\}$.
\item  $\partial(X)$ is good in $G$ if and only if $V(e_1)\subseteq X\cap U$, $V(e_2)\subseteq \overline{X}\cap W$, $|E[X\cap U,\overline{X}\cap W]|=2$ and $E[X\cap W,\overline{X}\cap U]$ contains only one edge (say $zw$, where $w\in X\cap W$). In addition, $zw$ is nonremovable in $G$, and both $G/\overline{X}$ and $G/X$ are near-bipartite graphs with removable doubleton $\{e_1,w\overline{x}\}$ and $\{e_2,zx\}$, respectively.
\end{enumerate}
\end{lemma}

\begin{proof}
Since $G$ is cubic and $\partial(X)$ is a 3-cut, $G/\overline{X}$ and $G/X$ are cubic, and $|\partial(X)|$ and $|X|$ have the same parity, so $|X|$ is odd. Since $G$ is 3-connected, so are $G/\overline{X}$ and $G/X$.
To show (i) and (iii), let $U'=X\cap U$, $W'=X\cap W$, $U''=\overline{X}\cap U$, $W''=\overline{X}\cap W$, $a=|U'|$ and $b=|W'|$.

(i) Suppose to the contrary that $a\geq b+2$. Since $G[U']$ has at most one edge $e_1$, $|\partial(X)|\geq 3a-2-3b\geq4$, a contradiction. (i) holds.

(iii) Suppose that $\partial(X)$ is good in $G$. Let $c=|E[U', \overline X]|$. Recall that $\partial(X)$ is a 3-cut. Then $c\leq 3$ and $|E[W', \overline X]|=3-c$.
By counting the number of edges in $E[U',W']$ in two ways, we have
$3a-2|E(U')|-c=3b-2|E(W')|-(3-c)$. Since  $a=b+1$, we have
\begin{equation}\label{41}c=|E(W')|-|E(U')|+3.\end{equation}
If $U'$ has at most one end of $e_1$, then $E(U')=\emptyset$. Since $c\leq 3$, we have $c=3$, $E(W')=\emptyset$ and $E[W',\overline{X}]=\emptyset$.
Thus $e_2\in E(W'')$. By (ii), $\partial(X)$ is tight in $G$, a contradiction.
So $e_1\in E(U')$. Analogously,  $e_2\in E(W'')$, so $E(W')=\emptyset$.
By (\ref{41}), we have $c=2$, so $|E[U',W'']|=2$ and $E[W',U'']$ contains only one edge, say $zw$, where $w\in W'$. Because every perfect matching of $G$ that contains $e_1$ contains $zw$, $zw$ is nonremovable in $G$.
Since both $G/\overline{X}-\{e_1,w\overline{x}\}$ and $G/X-\{e_2,zx\}$ are bipartite, Lemma \ref{lem:near-bipartite} implies that both $G/\overline{X}$ and $G/X$ are near-bipartite graphs with removable double $\{e_1,w\overline{x}\}$ and $\{e_2,zx\}$, respectively.

Conversely,  by Lemma \ref{lem:3-cut},  $\partial(X)$ is a separating cut of $G$. From (ii), we see that $\partial(X)$ is not tight in $G$, so it is good in $G$. (iii) holds.
\end{proof}

Let $G$ be a 3-connected  cubic graph. If $n=4$, then $G=K_4$. If $n\geq6$  and $G$ has a   nontrivial  3-cut $C$, then each $C$-contraction of $G$ is also a 3-connected  cubic graph that has strictly fewer vertices than $G$.  If either of the $C$-contractions has a nontrivial 3-cut, then the graph can be further decomposed into even smaller graphs. This procedure can be repeated until  we obtain a list of 3-connected  cubic graphs each of which is free of nontrivial 3-cuts. We refer to  it as a \emph{3-cut-decomposition} of $G$.
A 3-cut-decomposition is referred to as  a \emph{$K_4$-decomposition} if it results in a list of $K_4$s. If a graph has a $K_4$-decomposition,  by Lemma \ref{lem:cubic-brick-splicing}, it is a cubic brick and can be obtained by sequentially splicing  cubic bricks.

\begin{lemma}\label{lem-decomp}
Let $G$ be a cubic near-bipartite brick with at least six vertices. Then $G$ has a $K_4$-decomposition if and only if $G$ is a  tri-ladder.
\end{lemma}

\begin{proof}
Use induction on $n$. Note that $\overline{C_6}$ is the only graph obtained by splicing two $K_4$s, and it is a near-bipartite tri-ladder. Thus the result holds when $n=6$. Now suppose that $n\geq8$.

If $G$ is a tri-ladder, then $G$ has two vertex-disjoint triangles. Let $Y$ be the set of the vertices of one triangle. Then $\partial(Y)$ is a nontrivial 3-cut of $G$ and $G/\overline{Y}=K_4$.
By Lemma \ref{lem:3-cut} and the fact that  $G$ is a brick, $\partial(Y)$ is good in $G$. Let $G'=G/Y$. Then $G'$ is  a tri-ladder, which is also a cubic brick. By Lemma \ref{lem43}(iii), $G'$ is  near-bipartite. By the induction hypothesis, $G'$ has a $K_4$-decomposition, so does $G$.

If $G$ has a $K_4$-decomposition, then $G$ has a nontrivial 3-cut, say $\partial(X)$. Let $G_1=G/\overline{X}$ and $G_2=G/X$. Then   $G_i$, $i=1,2$,  has a $K_4$-decomposition, so  it is a cubic brick.
As above,  $\partial(X)$ is good.
Lemma \ref{lem43}(iii) implies that both $G_1$ and $G_2$ are  near-bipartite, and $\overline{x}$ and $x$ are incident with  edges in  removable doubletons of $G_1$ and $G_2$, respectively.
If $G_i$,  $i=1$ or $2$,  has at least six vertices, by the induction hypothesis, it is a tri-ladder. Then the contracted vertex, $\overline{x}$ or $x$, lies in a triangle of  $G_i$. If one of $G_1$ and $G_2$ is $K_4$, say $G_1$, then $G_2$ has at least six vertices, so it is a tri-ladder. Consequently, $G$ is a tri-ladder, which is a splicing of $G_2$ and $K_4$ at vertices in triangles.  If each of $G_1$ and $G_2$ has at least six vertices, then both of them are tri-ladders. Since a splicing of two tri-ladders at vertices in triangles is still a tri-ladder, $G$ is a tri-ladder. The result holds.
\end{proof}

\begin{lemma}\label{lem46}
Every near-bipartite tri-ladder  has exactly $\frac{n-6}{2}$ removable edges.
\end{lemma}

\begin{proof}
Let $G$ be a near-bipartite tri-ladder with $r^*$ rungs.
If $n=6$, then  $G=\overline{C_6}$, which has no removable edges. The result holds. Now assume that  $n\geq8$. Then $r^*\geq 1$.
By Theorem \ref{th:main1}, $G$ has at least $\frac{n-6}{2}$ removable edges. Note that $2r^*=n-6$. From the following claim, $G$ has exactly $\frac{n-6}{2}$ removable edges. We are done.

{\bf Claim.} Every edge $e$ that is not a rung is nonremovable in $G$.

To show this claim, we present some notions and notations.  Let $T$ and $T^*$ be the two disjoint triangles of $G$ with vertex set $\{u_0,v_0,w_0\}$ and $\{u,v,w\}$, respectively.
Assume that $u$ and $u_0$, $v$ and $v_0$, and $w$ and $w_0$ are ends of three ridges of $G$, which are referred to as $U$-ridge, $V$-ridge and $W$-ridge, respectively.
Let $f_i$ be the rung of $G$ with rank $i$, $1\leq i\leq r^*$.  We denote by  $u_i$, $v_i$, or $w_i$ the end of $f_i$ that lies in  $U$-ridge,  $V$-ridge, or $W$-ridge,  respectively.
For convenience, we refer to the edges in $T$ as rungs with rank zero.
By $u_iu_j$-subpath we mean the subpath of the $U$-ridge connecting $u_i$ and $u_j$. We use similar terminology for subpaths in  $V$-ridge and $W$-ridge.
Since $G$ is near-bipartite, $G$ has two edges  $e_1$ and $e_2$ such that $G-\{e_1,e_2\}$ is bipartite. Thus  one of $e_1$ and $e_2$ lies in $T$ and the other lies in $T^*$. Let $H=G-\{e_1,e_2\}$. Colour the vertices in the two color classes of $H$  white and black, respectively. We are now ready to prove the claim. We distinguish  two cases according to whether $e$ lies in triangles or not.

{\it Case 1.}  $e\in E(T\cup T^*)$. Suppose, without loss of generality, that  $e=uv$.
Assume that $r, s$ and $t$ are the maximum rank of rungs that has no end in $W$-ridge,  $U$-ridge and $V$-ridge, respectively. Thus these three rungs are  $u_rv_r,v_sw_s$ and $u_tw_t$.  For $i=r,s,t$,  let $d_i=r^*-i$. Then  $d_i\geq 0$.
Assume, without loss of generality, that $s\leq t$. Then we consider the neighbour of $v$ in the $V$-ridge, say $y$. Thus $y$ is either $v_r$ or $v_s$. Assume, without loss of generality, that $y$ is white. Then $v$ is black.
If $u$ is black, then $uv$ is either $e_1$ or $e_2$, so $uv$ is nonremovable in $G$. Suppose now that $u$ is white.  Recall that $r^{*}\geq1$.

First consider the case $y=v_r$ ($r\geq 0$). Recall that $0\leq d_r\leq r^{*}$.
Let $U'$ and $W'$ be the sets of all the black vertices in the $uu_r$-subpath and  the $ww_{r'}$-subpath, respectively,  where  $r'=r+1$ when $0<d_r<r^*$ and $r'=0$ when $d_r=r^*$ (i.e., $r=0$). If $d_r=0$, let $U'=\{u_r\}$ and $W'=\{w\}$. Let $B'=U'\cup W'\cup \{v_r\}$.
If $r\geq1$, since $u_rv_r\in E(H)$, $u_r$ is black.
If $r=0$ and $d_r$ is even,  since the $U$-ridge is a path of $H$, $u_0$ is black. In both cases,  $u_r\in B'$. Let $e'=u_rv_r$.
If $r=0$ and $d_r$ is odd, $u_0$ is white but $w_0$ is black, so $w_0\in B'$. Let $e'=w_0v_0$.
Consequently, $B'$ is a barrier of $G-uv$, which contains $V(e')$. By Lemma \ref{lem-barrier1}(ii), $e'$ is nonadmissible in $G-uv$, so $uv$ is nonremovable in $G$.

We now consider the case $y=v_s$ and $s\geq 1$. Then $w_s$ is black and $s\neq t$. Since $s\leq t$, we have $r^*= t>s$, so  $d_s\geq 1$. Then $u_{s+1}$ is black.
Therefore $v_s$ and all the black vertices  in the $ww_s$-subpath and  the $uu_{s+1}$-subpath form a barrier  of $G-uv$, which contains $v_s$ and $w_s$. As above, $uv$ is nonremovable in $G$.

{\it Case 2.}  $e\notin E(T\cup T^*)$. Then $e$ lies in a ridge, say the $V$-ridge. Let $v_{r^*+1}=v$.  Suppose, without loss of generality, that $e=v_jv_k$ ($0\leq j<k\leq r^*+1$), $v_j$ is white,  and the other end of the rung $f_j$ lies in the $U$-ridge. Then $f_j=u_jv_j$.
If $j=0$, then  $u_0w_0$ is nonadmissible in $G-v_0v_k$, so $v_0v_k$ is nonremovable in $G$. Suppose now that $j\geq 1$. Then $u_j$ is black.
Let $w_i$ be the vertex in the $W$-ridge with $i<j$ and $i$ as large as possible.
Then $f_{i+1}=u_{i+1}v_{i+1}$.

Suppose that $w_i$ is white.  When $j-i$ is odd,  $u_{i+1}$ is black and   $v_{i+1}$ is white.
Then $v_i$ is black, because $f_i=w_iv_i$ when $i\geq 1$, and $v_iv_{i+1}\in E(H)$ otherwise.
Let $B'$ be the set of the vertex $w_i$ and all the black vertices in the $u_{i+1}u_j$-subpath and the  $v_iv_j$-subpath.
Analogously, when $j-i$ is even, $u_{i+1}$ is white and  $v_{i+1}$ is black.
Then $u_i$ is black, because $f_i=w_iu_i$ when $i\geq 1$, and $u_iu_{i+1}\in E(H)$ otherwise.
Let $B'$ be the set of the vertex $w_i$ and all the black vertices in the $u_iu_j$-subpath and the $v_{i+1}v_j$-subpath. In both cases, $B'$ is a barrier of $G-v_jv_k$ containing two ends of  $w_iv_i$ or $w_iu_i$, so $v_jv_k$ is nonremovable in $G$.

Suppose that $w_i$ is black. Let $X$ be the set of all vertices in the $u_0u_j$-subpath, the $v_0v_j$-subpath and the $w_0w_i$-subpath.
Then $\partial(X)$  is a nontrivial 3-cut, which consists of the edge $v_jv_k$ and two edges incident with $u_j$ and $w_i$, respectively.
Since $G$ is a cubic brick, by Lemma \ref{lem:3-cut},  $\partial(X)$ is a separating cut but not tight, so it is good in $G$.
Recall that $v_j$ is white, and $u_j$ and $w_i$ are black.
Let $x_w$ and $x_b$ be the numbers of white and black vertices in  $X$, respectively, and $x_{wb}$ be the number of edges of $G$ each of which  connects a  white vertex and a black vertex in $X$.  By counting in two ways, we have $3x_{b}-2\geq x_{wb}\geq 3x_{w}-1-2$, so $x_{b}\geq x_{w}$.
By Lemma \ref{lem43}(iii),  $v_jv_k$ is nonremovable in $G$. The result holds.
\end{proof}

Next, we show that if $G$ is a near-bipartite brick other than $K_4$ that has exactly $\frac{n-6}{2}$ removable edges, then $G$ is a tri-ladder. Once this is proved. Theorem \ref{th:main2} holds. To do this, we need the following four lemmas.

\begin{lemma}[\cite{KCLL20}]\label{lem:essential-cubic-removable-doubleton}
In an essentially 4-edge-connected cubic brick, each edge is either removable
or lies in a removable doubleton.
\end{lemma}

\begin{lemma}[\cite{LFL20}]\label{lem:essentially-cubic-doubleton}
Suppose that $\mathscr{E}$ is a set of removable doubletons of an essentially 4-edge-connected cubic brick $G$ and $|\mathscr{E}|\geq2$. Then $G$ can be decomposed into balanced bipartite
vertex-induced subgraphs $G_i$ $(i=1,2,\ldots,|\mathscr{E}|)$ satisfying $E_G[V(G_j),V(G_k)]$
is a removable doubleton of $G$  if $|j-k|\equiv 1$ $(mod\quad\!\!\!\!|\mathscr{E}|)$, and  $E_G[V(G_j),V(G_k)]=\emptyset$ otherwise.
\end{lemma}

Here, a bipartite graph $G[A,B]$ is  balanced if $|A|=|B|$.

\begin{lemma}[\cite{LFL20}]\label{lem:two-adjacent-doubleton}
Suppose that $\{e_1,e_1'\}$ and $\{e_2,e_2'\}$ are removable doubletons of a cubic brick $G$.
If both $e_1$ and $e_2$ are incident with $v_0$, then $e_1'$ and $e_2'$ are adjacent, and $v_0u_0\in E(G)$, where $u_0$ is the common end of $e_1'$ and $e_2'$.
\end{lemma}

\begin{lemma}\label{prop-key}
Let  $\partial (X)$ be a  good 3-cut of a  cubic matching covered graph $G$.
Assume that $H=G/X\rightarrow x$ is an essentially 4-edge-connected cubic near-bipartite  brick other than $K_4$, and $\{xu,yv\}$ is a removable doubleton of $H$.
Let $(A,B)$ be the bipartition of $H-\{xu,yv\}$ such that $x,u\in A$.
Then each vertex of $B\setminus\{y,v\}$ is incident with at least two removable edges of $G$.
\end{lemma}

\begin{proof}
If $\{xu,yv\}$ is the only  removable doubleton of $H$, by Lemma \ref{lem:essential-cubic-removable-doubleton}, each edge of $E(H)\backslash\{xu,yv\}$ is removable in $H$. By Lemma \ref{lem41}, the result holds.
Now assume that $H$ contains $s$  removable doubletons, where $s\geq2$.

By Lemma \ref{lem:essentially-cubic-doubleton}, $H$ can be  decomposed into balanced bipartite subgraphs $H_1,H_2,\ldots,H_s$.
Let $A_i=V(H_i)\cap A$ and $B_i=V(H_i)\cap B$, $i=1,2\ldots,s$.
Let $x_1=x$, $u_s=u$, $y_1=y$, and $v_s=v$.
Moreover, we may suppose that $\{u_iy_{i+1},v_ix_{i+1}\}$, $i=1,2\ldots,s-1$, are $s-1$ removable doubletons of $H$ other than $\{x_1u_s,y_1v_s\}$, where $\{u_i,x_i\}\subseteq A_i$ and $\{y_i,v_i\}\subseteq B_i$.
Then $E_H[V(H_i),V(H_{i+1})]=\{u_iy_{i+1},v_ix_{i+1}\}$, $i=1,2\ldots,s-1$.
Let $w_1$ and $w_2$ be the two neighbors of $x_1$ other than $u_s$ in $H$.
Since $H$ is a brick, $H-\{w_1,w_2\}$ has a perfect matching, say $N_1$.
Note that $x_1u_s\in N_1$ because $d_H(x_1)=3$.
Since $\partial_G(X)$ is a  good 3-cut of  $G$, there exists a perfect matching of $G$, say $M$, containing each edge of $\partial_G(X)$.
Let $M'=(N_1\backslash\{x_1u_s\}\cup (M\setminus E(G[\overline{X}])$.
Then $M'$ is also a perfect matching of $G$.
Note that $w_1,w_2\in B_1\cup B_2$ and at most one of $w_1$ and $w_2$ lies in $B_2$.
We now consider the following two cases.

{\bf Case 1.} $w_1,w_2\in B_1$.
Since $|A_1\setminus\{ x_1\}|=|B_1\setminus\{w_1,w_2\}|+1$ and $E_H[A_1\setminus\{x_1\}, V(H)\setminus V(H_1)]=\{u_1y_2 \}$, we have $u_1y_2\in N_1$ and $v_1x_2\notin N_1$. So $u_1\neq x_1$.
Similarly, if $|s|\geq3$, we have $u_ty_{t+1}\in N_1$ and $v_tx_{t+1}\notin N_1$, $t=2,3,\ldots,s-1$.
Then $M'$ is a perfect matching of $G-(\cup_{i=1 }^{ s-1}\{v_ix_{i+1}\})$ and
$\cup_{i=1 }^{ s-1}\{u_iy_{i+1}\}\subseteq M'$.
Note that $G/\overline{X}$ is matching covered because $\partial_G(X)$ is a good cut of $G$.
For $i\in \{1,2,\ldots,s-1\}$, since $\{u_iy_{i+1},v_ix_{i+1}\}$ is a removable doubleton of $H$, $H-\{u_iy_{i+1},v_ix_{i+1}\}$ is matching covered. By Lemma \ref{lem:mc-splicing},  $G-\{u_iy_{i+1},v_ix_{i+1}\}$ is matching covered.
Recall that $M'$ is a perfect matching of $G-v_ix_{i+1}$ and $u_iy_{i+1}\in M'$.
Therefore, $G-v_ix_{i+1}$ is matching covered.
Then each edge of $\cup_{i=1 }^{ s-1}\{v_ix_{i+1}\}$ is removable in $G$.
By Lemma \ref{lem:essential-cubic-removable-doubleton}, each edge of $E(H)\setminus(\{x_1w_1,x_1w_2,x_1u_s, y_1v_s \}\cup (\cup_{i=1 }^{ s-1}\{u_iy_{i+1},v_ix_{i+1}\}))$ is removable in $H$, which is also removable in $G$  by Lemma \ref{lem41}.
This implies that each vertex of $B\setminus\{y_1,v_s\}$ is incident with at least two removable edges of $G$.
The result holds.

{\bf Case 2.} One of $w_1$ and $w_2$ lies in $B_2$, say $w_2$.
Then $w_1\in B_1$ and $x_1w_2=u_1y_2$. So $x_1=u_1$ and $w_2=y_2$. Since $\{x_1u_s,y_1v_s\}$ and
$\{x_1w_2, v_1x_2\}$ are removable doubletons of $H$, by Lemma \ref{lem:two-adjacent-doubleton}, we have $v_1=y_1$ and $x_1y_1\in E(H)$.
So $w_1=y_1=v_1$.
It follows that $|V(H_1)|=2$ because $E_H[V(H_1),V(H_2)]=\{x_1w_2,y_1x_2\}$ and $H$ is a cubic brick.
If $|s|=2$,
by Lemma \ref{lem:essential-cubic-removable-doubleton}, each edge of $E(H)\setminus\{x_1w_1,x_1w_2,x_1u_2,y_1v_2,y_1x_2\}$ is removable in $H$, which is also removable in $G$  by Lemma \ref{lem41}.
So the result holds.
We now consider $|s|\geq3$.
Since $|A_2|=|B_2\setminus\{w_2\}|+1$ and
$E_H[A_2, V(H)\setminus V(H_2)]=\{x_2y_1,u_2y_3\}$, we have $u_2y_3\in N_1$ and $v_2x_3\notin N_1$.
Similarly, if $|s|\geq4$, we have $u_ty_{t+1}\in N_1$ and $v_tx_{t+1}\notin N_1$, $t=3,4,\ldots,s-1$.
By the same reason as the above case, we can show that
each edge of $\cup_{i=2}^{ s-1}\{v_ix_{i+1}\}$ is removable in $G$, and each edge of $E(H)\setminus(\{x_1w_1,x_1w_2,x_1u_s, y_1v_s,y_1x_2\}\cup (\cup_{i=2}^{ s-1}\{u_iy_{i+1},v_ix_{i+1}\}))$ is removable in $G$.
So the result holds.
\end{proof}

\begin{lemma}\label{lem47}
If $G$ is a near-bipartite brick other than $K_4$ that has exactly $\frac{n-6}{2}$ removable edges, then $G$ is a  tri-ladder.
\end{lemma}

\begin{proof}
Let $\{e_1,e_2\}$ be a removable doubleton of $G$.
Since $G$ is a brick, $G$ is 3-connected. Thus, $\delta(G)\geq 3$.
By Theorem \ref{th:main1}, every vertex of $G$ is incident with at most two nonremovable edges and so at least one removable edge, except at most six vertices of degree three contained in two disjoint triangles of $G$. From the proof of Theorem \ref{th:main1} (in the first part of Section 3), we see that these two triangles contain $e_1$ and $e_2$, respectively.
Since $G$ has exactly $\frac{n-6}{2}$ removable edges, every vertex of $G$ other than the six vertices of degree three is incident with exactly one removable edge and two nonremovable edges.
Thus, $G$ is cubic. Note that $n\geq6$.
If $G$ has a $K_4$-decomposition, by Lemma \ref{lem-decomp}, $G$ is a tri-ladder.
We will show that $G$ has a $K_4$-decomposition in the following.

Since $G$ has triangles and $G\neq K_4$, it has nontrivial 3-cuts. Let $G^{*}$ be any graph obtained by a 3-cut-decomposition of $G$. Since $G$ is a cubic brick, $G^{*}$ is 3-connected and cubic. Thus, $G^{*}$ is an essentially 4-edge-connected cubic graph.
Let $\partial(X_1)$ be a nontrivial 3-cut of $G$ such that $x_1$ is a contraction vertex of $G^{*}$.
Since $G$ is a brick, $\partial(X_1)$ is good in $G$. By Lemma \ref{lem43}(iii), only one of $e_1$ and $e_2$ has both ends in $X_1$. This implies that $G^{*}$ has at most two contraction vertices.
By  Lemma \ref{lem:4-edge}, $G^{*}$ is a brick or a brace.
For the latter case, $|V(G^{*})|\geq6$ because $G^{*}$ is 3-connected and cubic. By Lemma \ref{lem:brace}, every edge of $G^{*}$ is removable in $G^{*}$. Since $G^{*}$ has at most two contraction vertices, there exists a vertex $v$ of $G^{*}$  that is incident with at most one contraction vertex of $G^{*}$.
By Lemma \ref{lem41}, $v$ is incident with at least two removable edges in $G$. This contradiction implies that $G^{*}$ is a brick.

By Lemma \ref{lem:cubic-brick-splicing}, each graph generated in the procedure of the 3-cut-decomposition of $G$ is a cubic brick.  Thus all the nontrivial 3-cuts used in this procedure are good.
By Lemma \ref{lem43}(iii), $G^{*}$ is a near-bipartite brick whose contraction vertices  are ends of the edges in its removable doubletons.
Then $G^{*}$ is an essentially 4-edge-connected cubic near-bipartite brick.
Recall that $x_1$ is a contraction vertex of $G^{*}$.
Let $\{x_1x_1',x_2x_2'\}$ be a removable doubleton of $G^{*}$. Suppose that $U'$ and $W'$ are two color classes of $G^{*}-x_1x_1'-x_2x_2'$ such that $x_1,x_1'\in U'$.
Suppose to the contrary that $G^{*}\neq K_4$. Then $|V(G^{*})|\geq6$.
If $G^{*}$ has exactly one contraction vertex $x_1$,
by Lemma \ref{prop-key}, each vertex of $W'\backslash \{x_2,x_2'\}$ is incident with at least two removable edges in $G$, a contradiction.
So $G^{*}$ has exactly two contraction vertices, one is $x_1$,  by Lemma \ref{lem43}(iii), the other is $x_2$ or $x_2'$, say $x_2$.
Let $G'=G/X_1$ and let $X_2$ be a subset of $V(G)$ such that $G^{*}=G'/X_2$. Then  $\partial(X_1)$ and $\partial(X_2)$  are good 3-cut in $G$ and  $G'$, respectively, and $G'$ is a cubic near-bipartite brick.
By Lemma \ref{prop-key}, each vertex of $U'\backslash\{x_1,x_1'\}$ is incident with at least two removable edges in $G'$, which are also removable in $G$ by Lemma \ref{lem41},  a contradiction. Thus $G^{*}=K_4$. The result holds.
\end{proof}

By Theorem \ref{th:main1}, Lemmas \ref{lem46} and \ref{lem47}, Theorem \ref{th:main2} holds.

\section*{Acknowledgements}
The authors would like to thank the anonymous referee for  helpful comments on improving the representation of the paper.
This work is supported by the National Natural Science Foundation of China (Nos. 12171440,   12271235   and 12371318), NSF of Fujian (Nos. 2021J01978, 2021J06029) and Fujian Alliance of Mathematics (No. 2023SXLMMS09).

\end{document}